\documentclass[12pt]{article}
\usepackage{amsthm,amsfonts,amsmath,amssymb}
\usepackage{amscd}

\usepackage{hyperref}

\newcommand{\eq}{\begin{equation}}
\newcommand{\en}{\end{equation}}

\renewcommand{\nabla}{{(-\Delta)}}

\newcommand{\giv}{\,|\,}

\newcommand{\prob}{\mathbb P}

\newcommand{\Nat}{\mathbb N}

\newcommand{\ed}{ \stackrel{d}{=}}

\newcommand{\PF}{{\cal V}}

\newcommand{\up}{\uparrow}
\newcommand{\down}{\downarrow}

\def\endpf{\hfill $\Box$ \vskip0.5cm}
\def \proof{\noindent{\it Proof.\ }}

\newtheorem{theorem}{Theorem}
\newtheorem{proposition}[theorem] {Proposition}

\newtheorem{corollary}[theorem]{Corollary}
\newtheorem{lemma}[theorem]{Lemma}

\theoremstyle{definition}
\newtheorem{definition}[theorem]{Definition}
\newtheorem{remark}[theorem]{Remark}
\newtheorem{example}[theorem]{Example}
\newtheorem{conjecture}[theorem]{Conjecture}
\newtheorem{problem}[theorem]{Problem}

\newcommand{\R}{\cal Z}

\begin{document}

\title{Coherent Permutations with Descent Statistic
and the Boundary Problem for the Graph of Zigzag Diagrams
}
\author{Alexander Gnedin\footnote{Utrecht University,
Mathematisch Instituut, PO Box 80010, 3508 TA Utrecht, The Netherlands,
gnedin@math.uu.nl\,.}\,\,\,\, and \setcounter{footnote}{6} Grigori
Olshanski\footnote{Institute for Information Transmission Problems, Bolshoy
Karetny 19,  Moscow 127994, Russia, olsh@online.ru\,.}}

\date{}

\maketitle

\begin{abstract}\noindent  The graph of zigzag diagrams is a close
relative of Young's lattice. The boundary problem for this graph amounts to
describing coherent random permutations with descent--set statistic, and is
also related to certain positive characters on the algebra of quasisymmetric
functions. We establish connections to some further relatives of Young's
lattice and solve the  boundary problem by reducing it to the classification of
spreadable total orders on integers, as recently obtained by Jacka and Warren.
\end{abstract}

\section{Introduction}

Young's lattice ${\cal Y}=\bigcup_{n\geq 0}\,{\cal Y}_n$ is the graded graph
whose vertices are Young diagrams, with the grading determined by the  number
of boxes in diagram, and with the neighbourship relation induced by the
inclusion of diagrams \cite[p. 288]{StanleyII}. We  write $\mu\nearrow\lambda$
when $\lambda$ is an immediate successor of $\mu$ in ${\cal Y}$, that is,
$\lambda$ is obtained from $\mu$ by adding a box. A path
$\lambda_0\nearrow\ldots\nearrow \lambda_n$ in $\cal Y$ starting with empty
diagram $\lambda_0=\varnothing$ encodes a standard Young tableau of the shape
$\lambda_n\in {\cal Y}_n$; such a path will be called a {\it standard\/} path
in $\cal Y$. The dimension function $d$ on $\cal Y$ which counts the number of
standard paths with a given end is uniquely determined by the recursion
\begin{equation}\label{recd}
d(\lambda)=\sum_{\mu:\,\mu\nearrow\lambda}d(\mu)\,,~~~~d(\varnothing)=1,
\end{equation}
whose explicit solution is given by Frobenius' formula or the famous hook
formula. Many questions of asymptotic combinatorics and representation theory
\cite{KerovDiss, KOO} involve the recursion
\begin{equation}\label{recY}
p(\mu)=\sum_{\lambda:\,\mu\nearrow\lambda}p(\lambda)\,,~~~~p(\varnothing)=1,
\end{equation}
which is dual to (\ref{recd}). The  boundary problem for $\cal
Y$ asks to describe the convex set of nonnegative solutions to
(\ref{recY}), especially the set $\partial{\cal Y}$ of extreme
solutions which we call the {\it boundary}.

\par The boundary problem has many facets.
One interpretation is that there is a bijection
$p\leftrightarrow\psi$ between nonnegative solutions $p$ to (\ref{recY}) and
linear functionals $\psi:{\rm Sym}\to{\mathbb R}$ on the algebra ${\rm Sym}$ of
symmetric functions, such that $\psi$ vanishes on a particular ideal and
assumes nonnegative values on the cone spanned by the basis of Schur functions.
By this bijection, the extreme solutions $p$ correspond exactly  to multiplicative functionals
$\psi$ (such $\psi$ will be called {\it characters.\/})  By another
interpretation, a nonnegative solution to (\ref{recY}) is a {\it probability
function\/}, which determines the distribution of a Markov chain $X=(X_n)$ with
$X_n\in {\cal Y}_n$, $n=0,1,2,\dots$, by  virtue of the formula

\begin{equation}\label{LMC}
{\mathbb P}(X_0=\lambda_0, \ldots, X_n=\lambda_n) =p(\lambda_n)\,,~~~~{\rm
for~~~} \varnothing=\lambda_0\nearrow\ldots\nearrow\lambda_n\in {\cal Y}_n\,.
\end{equation}

\par A central result on the boundary problem for Young's lattice is the
identification of $\partial {\cal Y}$ with the infinite `bi--simplex'
\begin{multline}\label{bisimplex}
\Delta^{(2)}=\{(\alpha,\beta): \alpha =(\alpha_j),\,\beta=(\beta_j),\\
\alpha_1\geq\alpha_2\geq\ldots\geq 0,~\beta_1\geq\beta_2\geq\ldots\geq 0,\\
\sum(\alpha_j+\beta_j)\leq 1\},
\end{multline}
as emerged in the work of Edrei on total positivity \cite{Edrei} and Thoma on
characters of the infinite symmetric group \cite{Thoma}, where the parameters
$(\alpha,\beta)$ appeared as the collections of poles and zeroes of the
generating function $\Sigma \,p(n)z^n$, with $p$ evaluated at one--row diagrams
$(n)$. See \cite{OkounkovDiss} for a related interpretation in terms of a
problem of moments. A probabilistic meaning of the parameters was discovered by
Vershik and Kerov \cite{VK-FA, VKSIAM}, who recognised in $\alpha_j$ and
$\beta_j$ the asymptotic frequencies of the $j$th largest row and the $j$th
largest column, respectively, for a diagram $X_n$ that follows the distribution
corresponding to $(\alpha,\beta)$. In Section \ref{Kerov} we exhibit yet
another construction of $\partial{\cal Y}$ based on the comultiplication in
$\rm Sym$, as  was communicated to the second author by Sergei Kerov many years
ago.

\par
The boundary problems for relatives of Young's lattice have been intensively
studied. Among them is the graph $\cal P$ of partitions of integers, which has
the same vertex set as $\cal Y$ but multiple edges, as dictated by the
branching of partitions of finite sets $[n]:=\{1,\ldots,n\}$; so that the
analogue of (\ref{recY}) involves some integer coefficients in the right--hand
side. From Kingman's work on partition structures \cite{Ki1, Ki2} we know that
the boundary $\partial{\cal P}$ is the simplex
$$\Delta=\{\alpha=(\alpha_j):  \alpha_1\geq\alpha_2\geq\ldots\geq 0,~
\sum \alpha_j\leq 1\}.$$ For $X=(X_n)$ a Markov chain on $\cal P$ corresponding
to $\alpha\in \Delta$, the parameter $\alpha_j$ has the same meaning of a row
frequency as in the case of Young's lattice, but the symmetry between rows and
columns, which holds for  ${\cal Y}$, breaks down for $\cal P$ and only the
longest column may have a positive frequency $1-\Sigma\,\alpha_j$. The boundary
problem for $\cal P$ may be also stated in terms of characters $\psi:{\rm
Sym}\to {\mathbb R}$ but this time $\psi$ must be nonnegative on a larger cone
in ${\rm Sym}$ spanned by the  monomial symmetric functions (whence the
embedding $\partial{\cal P}\hookrightarrow
\partial{\cal Y}$ which sends $\alpha\in \Delta$ to
$(\alpha,0)\in \Delta^{(2)}$). See \cite{KOO} for a  parametric
family of graded graphs which bridge between ${\cal Y}$ and
${\cal P}$.

\par
The graph of compositions ${\cal C}$ has compositions of integers as vertices,
and its standard paths of length $n$ correspond to ordered partitions of $[n]$.
This graph extends $\cal P$ in the sense that there is a projection ${\cal
C}\to {\cal P}$ which amounts to discarding the order of parts in composition.
The boundary problem for $\cal C$ is related to the characters $\psi:{\rm
QSym}\to{\mathbb R}$ of the algebra $ \rm QSym$ of quasisymmetric functions,
with the property that $\psi$ must be nonnegative on the cone spanned by the
basis of monomial quasisymmetric functions. In \cite{RCS} the boundary
$\partial{\cal C}$ was identified with the space $\cal U$ of open subsets
$U\subset [0,1]$ of the unit interval. The connection between the graphs
entails a projection $\partial{\cal C}={\cal U}\to\partial{\cal P}=\Delta$
which assigns to a generic open set $ U$ the decreasing sequence $\alpha$ of
sizes of its interval components; thus $U$ can be viewed as an arrangement of
Kingman's frequencies $\alpha_j$ in some order. There is an explicit {\it
paintbox} construction of Markov chains  on $\cal P$ and $\cal C$ by means of a
simple sampling scheme \cite{Ki1, Ki2, RCS, GP, CSP}.

\par
Thus, Kingman's boundary $\partial{\cal P}$ is a part of
the larger Edrei--Thoma--Vershik--Kerov  boundary $\partial{\cal
Y}$, and it is also  a shadow of the boundary $\partial{\cal C}$ of
the graph of compositions. This suggests to seek for a
combinatorial structure to complete the diagram
$$
\begin{CD}
\partial\,{\rm ?} @<<< \partial{\cal C}\\
@VVV @VVV\\
\partial{\cal Y} @<<<~ \partial{\cal P}
\end{CD}
$$
In this paper we argue that a good candidate to fill the gap is the {\it graph
of zigzag diagrams\/} $\cal Z$, because $\cal Z$ links to $\cal Y$ in a way
very similar to the relation between $\cal C$ and $\cal P$. The   vertices in
$\cal Z$ are again compositions, which we represent as zigzag diagrams, and the
configuration of edges is selected so that each standard path of length $n$
encodes a permutation of $[n]$ with a given set of descents. Alternatively, the
vertices of $\cal Z$ can be described as the words in a two--letter alphabet
with a natural interpretation of the relation $\mu\nearrow\lambda$, as explained in
the next section. In this interpretation, the graph $\cal Z$ is
known under the name of {\it subword order\/} and it was studied long ago, see
\cite{Viennot1983, Bjorner}. Together with the infinite binary tree, the graph
$\cal Z$ forms a dual pair of graphs in the sense of \cite{Fomin1994,
Fomin1995} (we thank Sergey Fomin for this remark). Recently, the graph $\cal
Z$ was considered in \cite{BS, Snellman}. However, all the papers listed above
are concerned with other problems, different from the boundary problem studied here. Notice
also that $\cal Z$ differs from other graphs of compositions \cite{PitmanPTRF,
KerovSub, Bosq} by the configuration of edges.

\par Although  $\cal Z$  cannot be literally projected on $\cal Y$, like $\cal
C$ on $\cal P$, their relation fits in the diagram of algebra homomorphisms
which respect  the cones spanned by the distinguished bases of symmetric
functions (s.f.) and quasisymmetric functions (q.s.f):

$$
\begin{CD}
({\rm QSym},~{\rm fundamental~q.s.f.})@>>> ({\rm QSym},~{\rm
monomial~q.s.f.})
\\
@AAA @AAA\\
({\rm Sym},~{\rm Schur~s.f.}) @>>>({\rm Sym},~{\rm
monomial~s.f.})
\end{CD}
$$

\par
We show that the comultiplication in $\rm QSym$ can be exploited to construct
the boundary $\partial{\cal Z}$, which turns out to be homeomorphic to the
space ${\cal U}^{(2)}$ of pairs $(U_\up,U_\down)$, where $U_\up$ and $U_\down$
are disjoint open subsets of $]0,1[$. The projection $\partial{\cal
Z}\to\partial{\cal Y}$ amounts to the separate record
$(U_\up,U_\down)\mapsto(\alpha,\beta)$ of sizes of interval components of the
open sets in decreasing order, while the embedding $\partial{\cal
C}\hookrightarrow\partial{\cal Z}$ is just $U\mapsto(U,\,\,\varnothing)$.
An extension of the paintbox scheme allows to describe the corresponding
Markov chains on $\cal Z$. In loose terms, the {\it oriented paintbox}
generates a random permutation by grouping $n$ integers in some number of
clusters and arranging the integers within each cluster in either increasing or
decreasing order. To show that our construction yields a complete description
of the boundary we reduce the boundary problem to a recent characterisation, due
to Jacka and Warren \cite{JW}, of
 {\it spreadable}  total orders on $\mathbb Z$
(a random total order is spreadable if its
probability distribution is invariant under increasing mappings
${\mathbb Z}\to{\mathbb Z}$).

\par The rest of the paper is organised as follows. In Section 2 we
introduce the graph of zigzag diagrams $\cal Z$. In Section 3 we explain how
$\cal Z$ is related to the algebra $\rm QSym$. In Section 4 we describe Kerov's
construction of the boundary $\partial{\cal Y}$ of the Young graph, and then in
Section 5 we apply Kerov's method to building a part of the boundary $\partial
{\cal Z}$, which we call the finitary skeleton of $\partial\cal Z$. In Section
6, using a natural continuity argument, we extend the finitary skeleton to a
larger set of boundary points, indexed by oriented paintboxes
$(U_\up,U_\down)$. In Section 7 we show that our construction actually yields
the whole boundary $\partial\cal Z$; this is the main result of the paper
(Theorem \ref{main}). The proof of the completeness claim amounts to a reduction to
the  remarkable result
by Jacka and Warren \cite{JW} on spreadable orders.
In Section 8 we formulate some complementary
results and open problems.

\section{The graph of zigzag diagrams}

Recall that a {\it composition\/} of $n$ is a finite sequence
$\lambda=(\lambda_1,\dots,\lambda_\ell)$ of strictly positive integers whose
sum $|\lambda|=\lambda_1+\dots+\lambda_\ell$ equals $n$.
Compositions have a useful
graphical representation, as in the following definition.

\begin{figure}[hbtp]

\unitlength 1mm

\linethickness{0.4pt}

\ifx\plotpoint\undefined\newsavebox{\plotpoint}\fi

\begin{picture}(102.02,28.87)(0,0)

\put(47.1,19.47){\framebox(4,4)[]{}}

\put(98.02,4.87){\framebox(4,4)[]{}}

\put(51.1,19.47){\framebox(4,4)[]{}}

\put(98.02,8.87){\framebox(4,4)[]{}}

\put(55.1,19.47){\framebox(4,4)[cc]{}}

\put(98.02,12.87){\framebox(4,4)[]{}}

\put(55.1,19.47){\framebox(0,.06)[]{}}

\put(98.02,12.87){\framebox(.06,0)[]{}}

\put(55.1,15.47){\framebox(4,4)[]{}}

\put(94.02,12.87){\framebox(4,4)[]{}}

\put(55.1,11.47){\framebox(4,4)[]{}}

\put(90.02,12.87){\framebox(4,4)[]{}}

\put(59.1,11.47){\framebox(4,4)[]{}}

\put(90.02,16.87){\framebox(4,4)[]{}}

\put(63.1,11.47){\framebox(4,4)[]{}}

\put(90.02,20.87){\framebox(4,4)[]{}}

\put(67.1,11.47){\framebox(4,4)[]{}}

\put(90.02,24.87){\framebox(4,4)[]{}}

\end{picture}

\caption{The zigzag shape  $\lambda=(3,1,4)$ and its conjugate
$\lambda'=(1^3,3,1^2)$}\label{rib}
\end{figure}
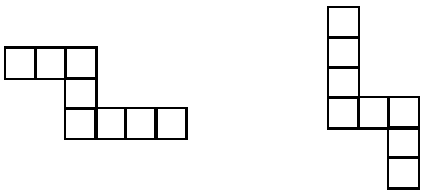

\begin{definition}\label{zigzag}
By a {\it zigzag diagram\/} (or simply {\it zigzag\/}) we understand a finite
collection of unit boxes such that the $j$th box is appended
either to the right or below  the $(j-1)$th box. (Thus, each zigzag comes with a canonical
enumeration of its boxes, from left to right and from top to bottom.) We will
identify compositions of $n=1,2\dots$ and zigzag diagrams with $n$ boxes: given
a composition $\lambda=(\lambda_1,\dots,\lambda_\ell)$, the corresponding
zigzag shape is drawn so that its $j$th row (counted from top
to bottom) consists of $\lambda_j$ boxes, $j=1,\dots,\ell$. Given a zigzag
diagram $\lambda$, the {\it conjugate\/} zigzag $\lambda'$ is obtained by
reflecting the zigzag $\lambda$ about the bisectrix of the first quadrant. See
Figure \ref{rib} for an example.
\end{definition}

\par The identification of compositions with zigzag diagrams is similar to the
identification of partitions with Young diagrams:
 both graphical
representations
unveil a symmetry that appears
by transposition of Young
diagrams or conjugation of zigzags.
Note also that the shape obtained by reflecting a zigzag $\lambda$ about the
horizontal axis is a ribbon Young diagram, as defined in \cite[p. 345]{StanleyII}.
\par
There exists an alternative encoding of compositions as {\it binary
words\/}. Take an alphabet with two symbols, say $\{+,-\}$, and assign to
zigzag $\lambda$ a word $w(\lambda)=w_1\dots w_{n-1}$ in this alphabet, as
follows: for $j=1,\dots, n-1$ let $w_j=+$ if the $j$th box in $\lambda$ is
separated from the $(j+1)$th box by a vertical edge, and let $w_j=-$ otherwise. To illustrate,
for $\lambda=(3,1,4)$, as on Figure 1, the word is $w(\lambda)=++--+++$. The
correspondence
 between zigzags with $n$ boxes
and binary words of length $n-1$ is a bijection (for $n=1$, the zigzag with one box is
represented by the empty word). As an obvious corollary we obtain that the total number
of zigzags with $n$ boxes equals $2^{n-1}$. In terms of binary words, the
conjugation $\lambda\to\lambda'$ amounts to switching $+\leftrightarrow-$ and
reading the resulting word in the reverse order.

\begin{definition}\label{word}
Let $\pi_n$ be a permutation of the set $[n]:=\{1,\dots,n\}$, written $\pi_n=\pi_n(1)\dots\pi_n(n)$
in the usual one--row notation.
The {\it zigzag shape\/} of
$\pi_n$, denoted  ${\rm zs}(\pi_n)$, is a zigzag $\lambda$ with $n$ boxes such
that in the corresponding binary word $w(\lambda)$, $w_j(\lambda)=+$ if
$\pi_n(j)<\pi_n (j+1)$ and $w_j(\lambda)=-$ if $\pi_n(j)>\pi_n(j+1)$ ($j=1,\dots,n-1$).
Thus, ${\rm zs}(\pi_n)=\lambda$ means that the numbers $\pi_n(1),\dots,\pi_n(n)$
inscribed consecutively in the boxes of $\lambda$ increase from left to right
along the rows and decrease down the columns. For example,
 ${\rm zs}(\pi_8)=(3,1,4)$ for $\pi_8=13842567$.
\end{definition}
\par
This definition can be rephrased as follows. Each permutation $\pi_n$ can be seen as a succession
 of increasing runs. The last position of each
increasing run, except the last run, is called a {\it descent\/}. For example,
the permutation $13842567$ breaks in three increasing runs $138, 4, 2567$ and  has two
descents $3$ and $4$ corresponding to the entries $\pi_8(3)=8$ and $\pi_8(4)=4$.
 The sizes of consecutive  runs comprise the composition
$\lambda={\rm zs}(\pi_n)$. Each of the descent positions in $\pi_n$
corresponds to
the rightmost box in some  row of $\lambda$
(different from the last row), or can be associated with
a horizontal edge separating two boxes.
Clearly, the zigzag shape of $\pi_n$ is  uniquely
 determined by the descent positions.
This encoding of descents
by zigzags goes back to MacMahon \cite[Ch. IV, Section
156]{MacMahon}.
\par
The zigzag shape of $\pi_n$ can also be determined from the inverse permutation $\pi^{-1}$
by recording  the sizes of increasing subsequences of consecutive integers.
For instance,  for $13842567$ the inverse permutation $\pi^{-1}=15246783$
breaks into subsequences $123,4,5678$ of sizes $(3,1,4)$.

\par
We need a few more definitions. For $\pi_n$ a permutation of $[n]$ and $\pi_m$
a permutation of $[m]$ with $m<n$, we say that the permutations are {\it
coherent\/} if the integers $1,\ldots,m$ appear in $\pi_n$ in the same relative
order as in $\pi_m$; we may also say that $\pi_n$ extends $\pi_m$, and that
$\pi_m$ is a restriction of $\pi_n$. In particular, $\pi_n$ extends $\pi_{n-1}$
if $\pi_n$ can be obtained by inserting element $n$ in the beginning, in the
end or between any two elements of the row $\pi_{n-1}(1)\ldots\pi_{n-1}(n-1)$.

\par
An infinite sequence of coherent permutations $\pi=(\pi_n)$ defines a total
order on $\Nat:=\{1,2,\dots\}$ and will be also called an {\it arrangement},
each $\pi_n$ is then the restriction of the arrangement $\pi$ to $[n]$. We will
use symbol $\lhd$ for arrangement as a relation on $\Nat$. Thus $i\lhd j$ means
that $\pi_n^{-1}(i)<\pi_n^{-1}(j)$ for every $n\geq \max(i,j)$. It is useful to
note that $\pi=(\pi_n)$ is uniquely determined by a single integer sequence
$(\pi_n^{-1}(n),\,n=1,2,\ldots)$ assuming values in $[1]\times[2]\times\cdots$.
We call $r_n:=\pi_n^{-1}(n)$ the {\it initial rank} of $n$, meaning that for
$\pi_n^{-1}(n)=r$ the integer $n$ is the $r$th $\lhd$--smallest in $[n]$.
Representing $\pi$ via initial ranks is rather convenient because passing from
$\pi_n$ to $\pi_{n+1}$ just amounts to appending $r_{n+1}$ to the row
$r_1,\ldots,r_n$.

\par Let ${\cal Z}_0=\{\varnothing\}$
be the one--element set which contains the empty zigzag
$\varnothing$, and for $n\geq 1$ let ${\cal Z}_n$ be the set of
all $2^{n-1}$ zigzags with $n$ boxes.

\begin{definition}\label{defZ}
The {\it graph of zigzag diagrams\/} has the set of vertices
${\cal Z}:=\cup_{n\geq 0} \,\,{\cal Z}_{n}$, and
the set of edges  defined by the
following neighbourship relation, further denoted $\nearrow$.
Let
$\varnothing\nearrow (1)$, and for $n\geq 2$ and two zigzags
 $\mu\in {\cal Z}_{n-1},~\lambda\in {\cal Z}_n$ let $\mu\nearrow\lambda$ if
$\mu={\rm zs}(\pi_{n-1})$ and $\lambda={\rm zs}(\pi_{n})$ for some coherent
permutations $\pi_{n-1}$ and $\pi_n$. We use the same symbol ${\cal Z}$ to
denote both the graph and its set of vertices.
\end{definition}

\par
\begin{figure}[hbtp]

\unitlength 1mm

\linethickness{0.4pt}

\ifx\plotpoint\undefined\newsavebox{\plotpoint}\fi

\begin{picture}(124.5,45.17)(0,0)

\put(70.25,29.5){\makebox(0,0)[cc]{(2)}}

\put(56.37,18.62){\makebox(0,0)[cc]{(3)}}

\put(75.84,18.8){\makebox(0,0)[cc]{(2,1)}}

\put(91.94,19.08){\makebox(0,0)[cc]{(1,2)}}

\put(51.87,7.12){\makebox(0,0)[cc]{(3,1)}}

\put(64.25,6.75){\makebox(0,0)[cc]{(1,3)}}

\put(75.39,7.42){\makebox(0,0)[cc]{$(2^2)$}}

\put(115.21,6.84){\makebox(0,0)[cc]{$(2,1^2)$}}

\put(124.5,7.09){\makebox(0,0)[cc]{$(1^4)$}}

\put(96.18,29.28){\makebox(0,0)[cc]{$(1^2)$}}

\put(108.83,18.8){\makebox(0,0)[cc]{$(1^3)$}}

\put(84.28,37.76){\makebox(0,0)[cc]{(1)}}

\put(89.89,7.26){\makebox(0,0)[cc]{$(1,2,1)$}}

\put(102.84,7.41){\makebox(0,0)[cc]{$(1^2,2)$}}

\put(39.62,7.12){\makebox(0,0)[cc]{(4)}}

\multiput(83,37)(-.06074286,-.0336){175}{\line(-1,0){.06074286}}

\multiput(85.5,36.87)(.0648125,-.0335625){160}{\line(1,0){.0648125}}

\multiput(92.37,18)(-.0335366,-.1036585){82}{\line(0,-1){.1036585}}

\multiput(92.12,18)(-.065952381,-.033730159){252}{\line(-1,0){.065952381}}

\multiput(92.62,17.75)(.04460581,-.03373444){241}{\line(1,0){.04460581}}

\multiput(64.37,9.12)(.105307692,.033653846){260}{\line(1,0){.105307692}}

\multiput(55,17.75)(-.052423077,-.033653846){260}{\line(-1,0){.052423077}}

\multiput(56.12,17.5)(-.03351852,-.07527778){108}{\line(0,-1){.07527778}}

\multiput(57.25,17.5)(.0337037,-.04100529){189}{\line(0,-1){.04100529}}

\multiput(57.87,17.37)(.07608491,-.03358491){212}{\line(1,0){.07608491}}

\multiput(76.12,16.5)(-.033333,-.208333){30}{\line(0,-1){.208333}}

\multiput(75,17.5)(-.04618257,-.03373444){241}{\line(-1,0){.04618257}}

\multiput(76.62,17.37)(.05538813,-.03365297){219}{\line(1,0){.05538813}}

\multiput(77.5,17.5)(.15298755,-.03373444){241}{\line(1,0){.15298755}}

\multiput(108.25,17.62)(-.03352941,-.05228758){153}{\line(0,-1){.05228758}}

\multiput(107.75,17.62)(-.07512605,-.03361345){238}{\line(-1,0){.07512605}}

\put(75.62,21.5){\line(0,1){0}}

\multiput(96.25,27.62)(-.03368794,-.04340426){141}{\line(0,-1){.04340426}}

\multiput(95.25,27.37)(-.11904192,-.03365269){167}{\line(-1,0){.11904192}}

\multiput(97.5,27.62)(.05829016,-.03367876){193}{\line(1,0){.05829016}}

\multiput(109.87,17.5)(.03369231,-.05961538){130}{\line(0,-1){.05961538}}

\multiput(111.37,17.75)(.05240664,-.03373444){241}{\line(1,0){.05240664}}

\put(83.87,43){\line(0,-1){3.12}}

\multiput(69.62,27.75)(-.0552093,-.03372093){215}{\line(-1,0){.0552093}}

\multiput(70.5,27.75)(.03355705,-.0385906){149}{\line(0,-1){.0385906}}

\multiput(71.5,28)(.10295337,-.03367876){193}{\line(1,0){.10295337}}

\put(83.88,45.17){\makebox(0,0)[cc]{$\varnothing$}}

\end{picture}

\caption{The beginning of the zigzag graph}\label{Rgraph}
\end{figure}
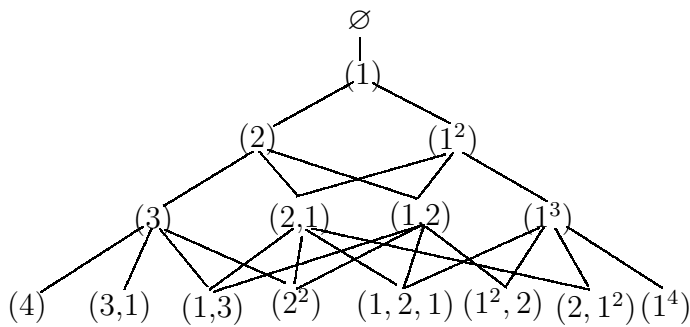

 Thus, $\cal Z$ is an infinite graded graph which also can be viewed as a graded
(or ranked) poset. Spelled out in detail, in terms of compositions, the
immediate followers of $(\mu_1,\ldots,\mu_\ell)$ are compositions
 $(\mu_1+1,\ldots,\mu_\ell)$,\dots,$(\mu_1,\ldots,\mu_\ell+1)$,
the composition $(1,\mu_1,\ldots,\mu_\ell)$, and all compositions
$$
(\mu_1,\ldots, \nu_j+1,\nu_{j+1},\ldots,\mu_\ell)
$$
obtained by splitting a
part $\mu_j$ in two positive parts $\nu_j+\nu_{j+1}=\mu_j$ and incrementing the
first of these two parts by $1$. Equivalently, in terms of zigzags, we split a
zigzag $\mu$ into two zigzag subdiagrams, $\mu(1)$ and $\mu(2)$, then according
as the edge between them is horizontal or vertical, we shift $\mu(2)$ by one
unit to the right or  down and finally add a new box in between. In the two
extreme cases when either $\mu(1)$ or $\mu(2)$ is empty, this means appending a
new box either above the first box of the zigzag or to the right of its last
box. Note that there are $n$ immediate followers $\lambda$ for any $\mu\in{\cal
Z}_{n-1}$. A few  levels of the graph are shown on Figure \ref{Rgraph}.

\begin{remark} In terms of binary words, the relation $\mu\nearrow\lambda$ (for
 $\mu\neq\varnothing$) holds if and only if  $w(\mu)$ is a subword of $w(\lambda)$, that
is, the word $w(\mu)$ can be obtained from the word $w(\lambda)$ by deleting
one of the symbols. It follows that for any zigzag $\lambda$ with
$|\lambda|\ge2$, the number of vertices $\mu\nearrow\lambda$ equals the total
number of plus-- and minus--clusters in $w(\lambda)$. In the `word'
realisation, the graph ${\cal Z}$
is known under
 the name  {\it subword
order\/}, see \cite{Viennot1983, Bjorner}.
\end{remark}

We define a {\it standard path of length $n$} in the graph $\cal Z$ as a
sequence $\lambda_0\nearrow\ldots\nearrow\lambda_n$ starting at the root
$\lambda_0=\varnothing$. Likewise, a standard infinite path is an infinite
sequence of neighbouring zigzags
$\lambda_0\nearrow\lambda_1\nearrow\lambda_2\ldots$.

\begin{proposition}\label{P1} For each $n$ and $\lambda_n\in {\cal Z}_n$
there is a canonical bijection between the permutations of $[n]$ with the
zigzag shape $\lambda_n$ and the standard paths in $\cal Z$ of length $n$
leading to $\lambda_n$. These bijections are consistent for various values of
$n$, hence define a bijection between the set of arrangements and the set of
infinite standard paths in $\cal Z$.
\end{proposition}

\proof Note that  the position of $n$ in a permutation of $[n]$
 is always a descent, unless this is the last position.
{F}rom this, inspecting all permutations $\pi_n$ extending a given permutation
$\pi_{n-1}$ with ${\rm zs}(\pi_{n-1})=\lambda_{n-1}$, we see that ${\rm
zs}(\pi_n)$ is determined by $\lambda_{n-1}$ and the position of $n$ in
$\pi_n$, and that the correspondence between the extensions of $\pi_{n-1}$ and
the successors of $\lambda_{n-1}$ in $\R$ is bijective. The claim follows by
induction in $n$.
\endpf

\noindent Proposition \ref{P1} becomes less evident when $\cal Z$ is introduced
as the subword order. In this form the statement
 was proved in  \cite{Viennot1983}.

\par The {\it dimension} function $d$ on ${\cal Z}$ counts the number of
permutations with a given descent set (equivalently, with a given zigzag
shape). It is easily seen that for a zigzag $\lambda\in{\cal Z}_n$, the
dimension $d(\lambda)$ equals the number of standard Young tableaux whose shape
is the ribbon Young diagram $\bar\lambda$, the image of $\lambda$ under the
reflection about a horizontal axis. Passing to the conjugate zigzag
$\lambda\mapsto\lambda'$ yields an involutive automorphism of $\cal Z$ hence
preserves the dimension.

\par
With the graph of zigzag diagrams we associate recursions of the same form
(\ref{recd}), (\ref{recY}) as for Young's lattice, but with understanding that
$\lambda,\mu\in {\cal Z}$. Obviously, we may equally  well write the
normalisation condition for (\ref{recY}) as $p((1))=1$.

\begin{definition}
Let $\PF$ denote the convex set of nonnegative solutions to recursion (\ref{recY})
where $\lambda$ and $\mu$ range over $\cal Z$.
We call the subset of extreme solutions the {\it boundary\/} of $\cal Z$ and
denote it $\partial{\cal Z}$.
\end{definition}

\par Note that $p(\lambda)\le1$ for any $\lambda\in\cal Z$. Using this observation
it is readily seen that $\PF$ is compact and metrisable in the product topology
(which corresponds to the pointwise convergence of functions on $\cal Z$). By a
well--known general theorem, it follows that $\partial{\cal Z}$ is a $G_\delta$
subset of $\PF$, hence a Borel space. Moreover, by the very definition $\PF$ is
a projective limit of finite--dimensional simplices. By
 \cite[p. 164,
Proposition 10.21]{Goodearl}
this implies  that the set $\PF$ is a Choquet simplex: that is,
any its point is {\it uniquely\/} representable as a convex mixture of extremes.
Thus, there is a canonical bijection between $\PF$ and
the set of all probability measures (`mixing
measures') on the space $\partial{\cal Z}$.

\par
There is  an equivalent probabilistic description of the set $\PF$ which is
useful for our purposes. Given a probability measure on the set of arrangements
we shall speak of a {\it random\/} arrangement $\Pi=(\Pi_n)$ which is a
sequence of coherent random permutations. By an easy application of
Kolmogorov's measure extension theorem, the law of $\Pi$ (which is  the {\it
joint} distribution of  the $\Pi_1,\Pi_2,\ldots$) is uniquely determined by the
{\it marginal} distributions of the $\Pi_n$'s. The following proposition
determines a family of random arrangements $\Pi=(\Pi_n)$ for which the descent
set of $\Pi_n$ is a sufficient statistic, for $n=1,2,\ldots$

\begin{proposition} The formula
\eq\label{Pp} {\mathbb P}(\Pi_n=\pi_n)=p({\rm zs}(\pi_n))\,, \qquad
n=1,2,\ldots \en establishes an affine homeomorphism $p\leftrightarrow\Pi$
 between  $\PF$ and the convex set of random arrangements
$\Pi=(\Pi_n)$ with the property that
for every $n$ all values of $\Pi_n$ with the same
zigzag shape are equally likely.
\end{proposition}

\proof This is a direct consequence of the definitions. Indeed, given
$p\in\PF$, we define, for each $n=1,2,\dots$ a
random permutation $\Pi_n$
 of $[n]$ with the probability law  $\prob_n(\Pi_n=\pi_n)=p({\rm zs}(\pi_n))$.
The relation (\ref{recY}) just means that, for every $n$, the measures
$\prob_n$ and $\prob_{n-1}$ are consistent with the projection
$\pi_n\mapsto\pi_{n-1}$ (removing $n$ from $\pi_n$), hence, by Kolmogorov's
theorem, $(\prob_n)$ uniquely determines a probability measure $\prob$ on the
set of arrangements, that is, a random arrangement $\Pi$. Clearly, such $\Pi$
satisfies (\ref{Pp}). Conversely, if a random arrangement $\Pi=(\Pi_n)$ with
distribution $\prob$ satisfies (\ref{Pp}) then (\ref{recY}) also holds by the
rule of addition for probabilities.
\endpf

\par
In slightly different terms,
not uncommon in this area of combinatorial probability \cite{CSP},
 we understand  each $p\in\PF$ as a {\it probability
function}, which by the formula
$$
{\mathbb P}(X_n=\lambda)=p(\lambda)d(\lambda)\,,~~~~~\lambda\in
{\cal Z}
$$
and an obvious analogue of (\ref{LMC})
determines the probability law of a transient
Markov chain $X=(X_n)$ whose
 state at time $n=0,1,\ldots$ is an element of ${\cal
Z}_n$.
The process $X$ viewed in the inverse time
$n=\ldots,2,1,0$ is again a Markov chain with transition
probabilities not depending on a particular choice of $p$, namely
$$
{\mathbb P}(X_{n-1}=\mu \mid X_n=\lambda)={\bf
1}(\mu\nearrow\lambda)\,d(\mu)/d(\lambda)\,,
$$
where ${\bf 1}(\cdots)$ equals $1$ when $\cdots$ is true and equals $0$
otherwise. By the Markov property, the probability law of $(X_n,\ldots,X_0)$ is
completely determined by the law of $X_n$, and the latter may be viewed as a
mixture of Dirac distributions on ${\cal Z}_n$.
\par Loosely speaking, a mixing measure on
$\partial{\cal Z}$ may be seen as a kind of initial distribution for the
time--reversed process $\ldots,X_2,X_1,X_0$ which  `starts at the infinite level of
$\cal Z$'.
By the general theory of Markov chains \cite{Kemeny}, $\partial{\cal
Z}$ is a subset of the larger {\it entrance Martin boundary\/} for the inverse chain,
and
in Section 8 we conjecture that the boundaries coincide.
Our purpose is to describe, as explicitly as possible, the boundary
$\partial{\cal Z}$
of the
graph of zigzag diagrams and the associated random processes on ${\cal Z}$.

\section{The algebra ${\rm QSym}$ }

The boundary problem for $\cal Z$ is intrinsically related to the {\it algebra
${\rm QSym}$ of quasisymmetric functions\/}. Usually ${\rm QSym}$ is introduced
 as a graded algebra of formal power series in infinitely many
indeterminates $x_1,x_2,\dots$. Specifically, elements of ${\rm QSym}$ are
finite linear combinations of the quasisymmetric monomial functions
$$
M_\lambda(x)=\sum_{1\le i_1<\dots<\lambda_\ell<\infty}x_{i_1}^{\lambda_1}\dots
x_{i_l}^{\lambda_\ell},
$$
where $\lambda=(\lambda_1,\dots,\lambda_\ell)$ is an arbitrary composition
(possibly, the empty one). The grading is defined by setting ${\rm
deg}M_\lambda=|\lambda|$.
Clearly, ${\rm QSym}$ contains the algebra ${\rm Sym}$ of symmetric functions
as a graded subalgebra.
\par
In ${\rm QSym}$ there is another homogeneous basis formed by the {\it
fundamental quasisymmetric functions\/} $F_\lambda$, see \cite{Gessel,
StanleyII, Thibon}:
$$
F_\lambda=M_\lambda+\sum_{\text{$\mu$ is thinner than $\lambda$}} M_\mu,
$$
where `$\mu$ is thinner than $\lambda$' means that the composition $\mu$ can be
obtained from $\lambda$ by splitting some  parts of the composition $\lambda$ into
smaller parts, so that at least one part $\lambda_j$ is replaced by a
nontrivial composition of $\lambda_j$. For example,
$$
F_{(2,2)}=M_{(2,2)}+M_{(1,1,2)}+M_{(2,1,1)}+M_{(1,1,1,1)}\,.
$$
Note that ${\rm deg}F_\lambda=|\lambda|$. The role of functions $F_\lambda$ is
similar to that of Schur functions in the algebra ${\rm Sym}$.

\par For most of our
purposes it suffices to simply regard ${\rm QSym}$ as a commutative unital
$\mathbb R$--algebra with distinguished basis $\{F_\lambda\,,\lambda\in{\cal Z}\}$ indexed
by compositions,   including the unity $F_\varnothing=1$.
The multiplication rule in this basis is the following.
For $\pi_k$ a permutation of
$[k]$ and $\pi_l$ a permutation of $[l]$ let ${\rm shuffle}(\pi_k,\pi_l)$ be
the set of permutations $\pi_{k+l}$ of $[k+l]$ such that $\pi_k$ and
$\pi_{k+l}$ are coherent, and the relative order of integers $k+1,\ldots,k+l$
in $\pi_{k+l}$ is the same as in the row $k+\pi_l(1),\ldots,k+\pi_l(l)$. Then
for $\mu={\rm zs}(\pi_k)$ and $\nu={\rm zs}( \pi_l)$  the  product is defined
as
\begin{equation}\label{mult}
F_{\mu}F_\nu=\sum_{\pi_{k+l}\,\in \,{\rm shuffle}(\pi_k,\pi_l)}
F_{{\rm zs}(\pi_{k+l})}.
\end{equation}
By \cite[p. 483, Exercise 7.93]{StanleyII} this  operation is commutative,
associative and indeed well defined. That is to say, the definition  is not sensitive to the
choice of permutations $\pi_k$, $\pi_l$ provided they have  given zigzag
shapes.  It is worth noting that the multiplication rule of the fundamental
quasisymmetric functions is much simpler than that of Schur functions.

\par
Of  special importance for us is the instance of (\ref{mult}) in
which $\pi_k=\pi_{n-1}$ is a permutation of $[n-1]$ and
$\pi_l=\pi_1$ is the sole permutation of $[1]$. Since $\pi_n\in
\,{\rm shuffle}(\pi_{n-1},\pi_1)$ just means that $\pi_n$
extends $\pi_{n-1}$, we have
\begin{equation}\label{mult1}
F_\mu F_{(1)}=\sum_{\lambda:\, \mu\nearrow\lambda} F_\lambda\,.
\end{equation}
\par
The graph $\cal Z$ is a {\it multiplicative graph\/} associated to ${\rm QSym}$
in the sense of Vershik--Kerov \cite{KerovDiss}, meaning that (\ref{mult1})
mimics the branching in $\cal Z$, and that the cone spanned by
$\{F_\lambda,~\lambda\in {\cal Z}\}$ is stable under multiplication. A
consequence of the multiplicativity is that the formula
\begin{equation}\label{psi-p}
\psi(F_\lambda)=p(\lambda)
\end{equation}
establishes a bijection
 between probability functions  and linear
functionals $\psi:{\rm QSym}\to {\mathbb R}$ with the properties
\begin{itemize}
\item[(i)] $\psi(1)=1\,$,
\item[(ii)] $\psi(F_\lambda)\geq 0\,$ for $\lambda\in {\cal Z\,}$,
\item[(iii)] $\psi(F_\lambda F_{(1)})=\psi(F_\lambda)\,$ for $\lambda\in {\cal Z}\,$.
\end{itemize}

Let ${\rm QSym}^*$ be the algebraic dual space to the vector space ${\rm
QSym}$. Denoting $K$ the cone generated by $\{F_\lambda,\,\lambda\in {\cal
Z}\}$, condition (ii) says that $\psi$ belongs to the dual cone $K^*\subset{\rm
QSym}^*$.

\begin{proposition}\label{P2}
Under the correspondence {\rm (\ref{psi-p})},   $p\in
\partial{\cal Z}$ if and only if $\psi$ is multiplicative, in
the sense that
\begin{itemize}
\item[\rm(iv)]  $\psi(FG)=\psi(F)\psi(G)$  for  $F,G\,\in {\rm QSym}\,.$
\end{itemize}
\end{proposition}

\proof  By virtue of the multiplicativity of $\cal Z$, this follows from  the
general `ring theorem' \cite{VK-Roman, KerovDiss}. See also Proposition
\ref{ring} in the Appendix below. We apply this result to the quotient of
${\rm QSym}$ by the ideal generated by $F_{(1)}-1$, with  the cone being the
image of the cone $K$ by the canonical projection.
\endpf

A functional $\psi$ which fulfills (i)--(iv) will be called a {\it character\/}
of ${\rm QSym}$. The set of characters is a subset of the space ${\rm QSym}^*$
endowed with the weak
topology.

\begin{corollary} The boundary $\partial{\cal Z}$ and the set of
characters of $\,{\rm QSym}$ are compact topological spaces homeomorphic to each other.
\end{corollary}

\proof Since the weak topology on ${\rm QSym}^*$ coincides with the topology of
simple convergence on the basis elements $F_\lambda$, the correspondence
(\ref{psi-p}) is a homeomorphism. On the other hand, by Proposition \ref{P2},
the boundary is closed in ${\rm QSym}^*$, hence in $\PF$. Since $\PF$ is
compact, so is the boundary.
\endpf

\par
Besides multiplication, there are two other important operations in ${\rm
QSym}$ of involution and  comultiplication, which we introduce next.

\begin{definition}\label{inv}
Let ${\rm inv}: {\rm QSym}\to{\rm QSym}$ be the involutive linear map defined
by ${\rm inv}(F_\lambda)=F_{\lambda'}$, where $\lambda'$ is the conjugate of
$\lambda\in {\cal Z}$.
\end{definition}

\par Using (\ref{mult}) it is readily checked that $\rm inv$ is an algebra
automorphism. Moreover, ${\rm inv}$ extends the canonical involution $\omega$
of the algebra ${\rm Sym}$. See Subsection 8.1 for more detail.

\begin{definition}\label{comult}
(i) The {\it comultiplication\/} is a graded algebra homomorphism $\delta:{\rm
QSym}\to{\rm QSym}\otimes{\rm QSym}$. In the basis of monomial functions
$\delta$ can be defined by splitting the family of variables $x$ into
two disjoint infinite subsets $x'$ and $x''$. A generic quasisymmetric
function $F$ becomes then a function $F(x',x'')$ of two families of variables
$x'$ and $x''$. Because $F$ is
separately quasisymmetric in both $x'$ and  $x''$, the function $F$ is
 an element of ${\rm QSym}\otimes{\rm QSym}$.

(ii) Equivalently, $\delta$ can be defined on the basis of fundamental
quasisymmetric functions as
\begin{equation}\label{delta}
\delta(F_\lambda)=\sum_{\mu,\nu: \,\,\lambda=\mu \sqcup\nu} F_\mu\otimes F_\nu\,,
\end{equation}
where, for $\lambda\in {\cal Z}_n$, the sum is over all $n+1$ splittings
$\lambda=\mu \sqcup\nu$ in two zigzags
$\mu$ and $\nu$ one of which can be empty.
\end{definition}

The fact that (\ref{delta}) is indeed an algebra homomorphism can be checked
directly using (\ref{mult}). Note also that (\ref{delta}) is dual to the rule
  \cite[Equation (67)]{Thibon}. From (i) it follows that $\delta$ extends the
comultiplication in ${\rm Sym}$, which can be defined in exactly the same way.
However, in contrast to the algebra ${\rm Sym}$, the comultiplication in ${\rm
QSym}$ is not cocommutative, since switching the factors in ${\rm QSym}\otimes{\rm
QSym}$ does not preserve $\delta$.

\begin{proposition}\label{projection}
The canonical embedding\/ ${\rm Sym}\hookrightarrow{\rm QSym}$ induces a
projection of the boundaries $\partial{\cal Z}\to\partial{\cal Y}$.
\end{proposition}

\proof Suppose $\psi:{\rm QSym}\to{\mathbb R}$ satisfies (i)--(iv). The Schur
functions $S_\lambda$, $\lambda\in {\cal Y}$, expand with nonnegative
coefficients in the basis of fundamental quasisymmetric functions \cite[Theorem
7.19.7]{StanleyII}, hence the cone generated by $\{S_\lambda,~\lambda\in {\cal
Y}\}$ is thinner than the cone $K={\rm conv}\{F_\lambda,~\lambda\in{\cal Z}\}$
, hence $\psi$ is nonnegative on the Schur functions. That $\psi$ is also a
character of ${\rm Sym}$ follows from Proposition \ref{P2}, since
 $F_{(1)}=S_{(1)}$ and the embedding is a homomorphism.
\endpf

\noindent Thus,  we can expect {\it a priori} that the parameters
$(\alpha,\beta)$ of $\partial{\cal Y}$ will enter, in some way, in the
description of the boundary of the graph of zigzag diagrams.

\section{Kerov's construction of $\partial{\cal Y}$}
\label{Kerov}

The algebra $\rm Sym$ of symmetric functions with  the basis $\{S_\lambda,
~\lambda\in {\cal Y}\}$ of Schur functions underlies Young's graph $\cal Y$,
which is multiplicative and has the branching  read from the relation
$$
S_\mu S_{(1)}=\sum_{\lambda:\,\mu\nearrow\lambda} S_\lambda\,.
$$
Replacing the functions $F_\lambda$
by the functions $S_\lambda$ in the formula
(\ref{psi-p}) and the subsequent discussion,
we identify the elements of the
boundary $\partial\cal Y$ with the characters $\psi$ of the algebra ${\rm
Sym}$.  There exists a simple construction of the characters which is based on
the comultiplication in the algebra ${\rm Sym}$. This construction, due to
Kerov, did not appear in the literature before. We present it in detail because
in the next sections it will be used as a prompt for a more involved
construction of the characters of ${\rm QSym}$.
\par As mentioned above, the comultiplication $\delta:{\rm Sym}\to {\rm
Sym}\otimes{\rm Sym}$ is consistent with the comultiplication in ${\rm QSym}$
under the embedding ${\rm Sym}\to{\rm QSym}$. We have
\begin{equation}\label{deltaS}
\delta(S_\lambda)=\sum_{|\mu|+|\nu|=|\lambda|}c(\lambda;\mu,\nu)S_\mu\otimes
S_\nu\,,
\end{equation}
where the coefficients $c(\lambda;\mu,\nu)$ coincide with the structure
constants of the multiplication, the so-called Littlewood--Richardson coefficients. This
coincidence is a manifestation of the fact that ${\rm Sym}$ is a self--dual
Hopf algebra, see \cite[Ch. I, section 5, Example
25]{Macdonald}. The fact that the coefficients $c(\lambda;\mu,\nu)$ are
nonnegative  plays a crucial role in the construction.

\par
\begin{proposition}\label{prop1}
Let $\psi_1$ and $\psi_2$ be two characters of\/ ${\rm Sym}$ and let $\omega_1$
and $\omega_2$ be two
nonnegative numbers whose sum equals $1$. Define a
linear functional $\psi:{\rm Sym}\to{\mathbb R}$ by
$$
\psi(S_\lambda)=\sum_{|\mu|+|\nu|=|\lambda|}c(\lambda;\mu,\nu)\omega_1^{|\mu|}
\omega_2^{|\nu|} \psi_1(S_\mu)\psi_2(S_\nu).
$$
Then $\psi$ is a character too.
\end{proposition}

\proof Let us check properties (i)--(iv) above (with $F$--functions replaced by
$S$--functions). Property (i) is trivial because $\delta(1)=1\otimes1$.
Property (ii) follows from the nonnegativity of coefficients
$c(\lambda;\mu,\nu)$. Regarding (iii), observe that
$$
\delta (S_{(1)})=S_{(1)}\otimes1+1\otimes S_{(1)}.
$$
Since $\psi_1(S_{(1)})=1$, $\psi_2(S_{(1)})=1$, and $\omega_1+\omega_2=1$, this
implies $\psi(S_{(1)})=1$. Hence (iii) follows from (iv). To prove (iv),
introduce a one--parameter family of algebra endomorphisms $r_\omega:{\rm
Sym}\to{\rm Sym}$,  where $\omega\in\mathbb R$, by setting
\begin{equation}\label{r}
r_\omega(S_{\lambda})=\omega^{|\lambda|}S_\lambda
\end{equation}
($r_\omega$ is an automorphism for $\omega\ne0$). By the very definition of
$\psi$, it coincides with the superposition of $\psi_1\otimes\psi_2$,
$r_{\omega_1}\otimes r_{\omega_2}$, and $\delta$. Therefore, $\psi$ is
multiplicative.
\endpf

Set $\delta^{(2)}=\delta$ and define $\delta^{(3)}:{\rm Sym}\to{\rm
Sym}^{\otimes 3}$ as the superposition map
$$
\delta^{(3)}=(\delta\otimes{\rm id})\circ\delta^{(2)}=({\rm
id}\otimes\delta)\circ\delta^{(2)}.
$$
Here ${\rm id}:{\rm Sym}\to{\rm Sym}$ is the identity map and the second
equality above follows from the associativity of the comultiplication.
Continuing this procedure we recursively define the iterated map
$\delta^{(k)}:{\rm Sym}\to{\rm Sym}^{\otimes k}$ for any $k=2,3,\dots$.
Clearly, this is a graded algebra morphism. Just as in the case $k=2$,
$\delta^{(k)}$ can be defined by splitting the family of variables into $k$
disjoint subsets.

Now we can generalise Proposition \ref{prop1}.

\begin{proposition}\label{prop2}
For any $k=2,3,\dots$ let $\psi_1, \dots,\psi_k$ be characters of\/ ${\rm Sym}$
and let $\omega_1, \dots, \omega_k$ be real nonnegative numbers whose sum
equals $1$. Define a linear functional $\psi:{\rm Sym}\to{\mathbb R}$ as the
superposition
\begin{equation}\label{M}
(\psi_1\otimes\dots\otimes\psi_k) \circ(r_{\omega_1}\otimes\dots\otimes
r_{\omega_k})\circ\delta^{(k)},
\end{equation}
where the endomorphisms $r_\omega$ were introduced in\/ {\rm(\ref{r})}. Then
$\psi$ is a character too.
\end{proposition}

\proof By the definition of $\delta^{(k)}$, all the coefficients in the
expansion of $\delta^{(k)}(S_\lambda)$ in the natural basis of ${\rm
Sym}^{\otimes k}$ (the tensor products of $S$--functions) are nonnegative. We
also observe that
$$
\delta^{(k)}(S_{(1)})=S_{(1)}\otimes1\otimes\dots\otimes1+1\otimes
S_{(1)}\otimes1\otimes\dots\otimes1+\dots+1\otimes\dots\otimes1\otimes S_{(1)}.
$$
Then we argue as in the proof of Proposition \ref{prop1}.
\endpf

\begin{definition}\label{Mmix}
The above operation will be called {\it M--mixing\/} (M for
`multiplicative'). We will say that the character (\ref{M}) is obtained by
M--mixing the characters $\psi_1,\dots,\psi_k$ in proportions
$\omega_1,\dots,\omega_k$.
\end{definition}

Observe that the characters of $\rm Sym$ are uniquely specified by their values
on the complete homogeneous symmetric functions $h_n=S_{(n)}$ which are
generators of the algebra. Using the generating series $H(t)=1+\sum_n h_nt^n$
with formal parameter $t$ we set for a character $\psi$
$$
\psi(H(t))=1+\sum_{n=1}^\infty \psi(h_n)t^n.
$$

\begin{proposition}\label{psi(H)} If $\psi$ results from  M--mixing the characters
$\psi_1,\dots,\psi_k$ in proportions $\omega_1,\dots,\omega_k$ then
$$
\psi(H(t))=\psi_1(H(\omega_1 t))\dots\psi_k(H(\omega_k t)).
$$
\end{proposition}

\proof Indeed, this follows from the definition of M--mixing and the
well--known formula
$$
\delta^{(k)}(H(t))=H(t)\otimes\dots\otimes H(t) \qquad\text{($k$ times)}.
$$
which follows, e.g., from \cite[Ch. I, Section 5, Example 25
(a)]{Macdonald}.
\endpf

The next proposition introduces two elementary characters
which will be used  to generate further characters
by M-mixing.

\begin{proposition}\label{elemchar}
 The linear functionals $\psi_+:{\rm Sym}\to\mathbb R$ and
$\psi_-:{\rm Sym}\to\mathbb R$ defined by
$$
\psi_+(S_\lambda)=\begin{cases} 1  &\text{if $\lambda=(n)$}\\ 0, & \text{if
$\lambda\ne(n)$}
\end{cases},
\qquad \psi_-(S_\lambda)=\begin{cases} 1 &\text{if $\lambda=(1^n)$}\\ 0, &
\text{if $\lambda\ne(1^n)$} \end{cases}
$$
are characters. Here $(n)$ is the one--row diagram and $(1^n)=(n)'$ is the
conjugate one--column diagram, $n=0,1,2,\dots$.
\end{proposition}

\proof The function $p$ on $\cal Y$ corresponding to $\psi_+$ is the
characteristic function of the set of one--row diagrams. Clearly, $p$ is an
extreme element in the convex set of solutions to (\ref{recY}). Therefore, by
the ring theorem (see Proposition \ref{ring}), $\psi_+$ is multiplicative. For
$\psi_-$ the argument is similar.

An alternative proof consists in checking directly that $\psi_+$ and $\psi_-$
are multiplicative (other properties being obvious). For $\psi_+$ this reduces
to the following two claims. Firstly, $h_mh_n=h_{m+n}+\dots$, where the remaining
term does not involve one--row $S$--functions (indeed, it is well known that
the remaining term is a sum of two--row functions). Secondly, if at least one of
the two diagrams $\mu$ and $\nu$ has 2 or more rows then the expansion of $S_\mu
S_\nu$ does not involve one--row $S$--functions (indeed, as is well known, if
$S_\lambda$ enters the expansion of $S_\mu S_\nu$ then $\lambda$ must contain
both $\mu$ and $\nu$). For $\psi_-$ one can apply a similar argument or use the
canonical involution of the algebra ${\rm Sym}$ which sends $\psi_+$ to
$\psi_-$.
\endpf

Now we proceed to the explicit construction of characters. Given an arbitrary
collection of real nonnegative numbers $\alpha_1,\dots,\alpha_k$ and
$\beta_1,\dots,\beta_l$ whose total sum equals 1, we take the M--mixing of $k$
copies of $\psi_+$ and $l$ copies of $\psi_-$ in proportions
$\alpha_1,\ldots,\alpha_k$ and $\beta_1,\ldots,\beta_l$, respectively. Because
the comultiplication in ${\rm Sym}$ is cocommutative, the order in which the
couples $(\alpha_i, \psi_+),\,(\beta_j,\psi_-)$ are taken plays  no role. This
is reflected in the symmetry of formula (\ref{psi(H(t))}) below with respect to separate
permutations of $\alpha_i$'s and $\beta_j$'s.

\begin{proposition} For the character $\psi$ obtained by the above construction
we have
\begin{equation}\label{psi(H(t))}
\psi(H(t))=\frac{\prod\limits_{j=1}^l(1+\beta_j t)}
{\prod\limits_{i=1}^k(1-\alpha_i t)}\,.
\end{equation}
\end{proposition}

\proof By the definition of $\psi_+$,
$$
\psi_+(H(t))=\frac1{1-t}.
$$
Likewise, denoting by $E(t)$ the generating series for the elementary symmetric
functions $e_n=S_{(1^n)}$ we have
$$
\psi_-(E(t))=\frac1{1-t}.
$$
Together with the classical relation $H(t)E(-t)=1$ this implies
$$
\psi_-(H(t))=1+t.
$$
Now the result follows from Proposition \ref{psi(H)}.
\endpf

\par Arranging the parameters in decreasing order we may identify
$(\alpha_1,\dots,\alpha_k$; $\beta_1,\dots,\beta_l)$ with the element
$$
(\alpha_1,\dots,\alpha_k, 0,0,\dots;\beta_1,\dots,\beta_l,0,0,\dots)
\in\Delta^{(2)}.
$$
Let $\Delta^{(2)}_{\rm fin}\subset\Delta^{(2)}$ stand for the subset of all
elements of such a kind, that is, $\Delta^{(2)}_{\rm fin}$ consists of all
couples $(\alpha,\beta)\in\Delta^{(2)}$ with finitely many nonzero coordinates
whose total sum equals 1. We have shown that to any
$(\alpha,\beta)\in\Delta^{(2)}_{\rm fin}$ there is a character
$\psi_{\alpha,\beta}$ such that
$$
\psi_{\alpha,\beta}(H(t))=\frac{\prod\limits_{j=1}^\infty(1+\beta_j t)}
{\prod\limits_{i=1}^\infty(1-\alpha_i t)}\,,
$$
where the products are actually finite. We call $\{\psi_{\alpha,\beta}:
(\alpha,\beta)\in\Delta^{(2)}_{\rm fin} \}$ the {\it finitary skeleton\/} of
the boundary $\partial\cal Y$.
\par
Let ${\rm Sym}^*$ be the algebraic dual space to the vector space ${\rm Sym}$.
We equip ${\rm Sym}^*$ with the weak topology. It is readily seen that the
characters form a closed subset of ${\rm Sym}^*$. Moreover, the induced
topology on the set of characters coincides with the product topology by the embedding
$\psi\mapsto (\psi(h_1),\psi(h_2),\dots)$.

\begin{proposition} The closure of the finitary skeleton in ${\rm Sym}^*$
consists of the characters which are indexed by arbitrary elements
$(\alpha,\beta)\in\Delta^{(2)}$ and are specified by
\begin{equation}\label{psi}
\psi_{\alpha,\beta}(H(t))=e^{\gamma
t}\,\frac{\prod\limits_{j=1}^\infty(1+\beta_j t)}
{\prod\limits_{i=1}^\infty(1-\alpha_i t)}\,, \qquad{\rm where~~~}
\gamma:=1-\sum_{i=1}^\infty\alpha_i-\sum_{j=1}^\infty\beta_j\,.
\end{equation}
\end{proposition}

\proof Let ${\mathbb R}^\infty$ be the product of countably many copies of
${\mathbb R}$ taken with the product topology. The set $\Delta^{(2)}$ may be
viewed as a subset of ${\mathbb R}^\infty\times{\mathbb R}^\infty$.  To any
$(\alpha,\beta)\in\Delta^{(2)}$ we assign a multiplicative functional
$\psi_{\alpha,\beta}:{\rm Sym}\to\mathbb R$ by making use of (\ref{psi}).
Observe that $\Delta^{(2)}$ is compact in the induced topology and
$\Delta^{(2)}_{\rm fin}$ is  dense in $\Delta^{(2)}$. Therefore, it suffices to
prove that the quantities $\psi_{\alpha,\beta}(h_n)$ are continuous functions
in $(\alpha,\beta)\in\Delta^{(2)}$. It is more convenient to consider the
values on another set of generators, the Newton power sums
$p_1,p_2,\ldots\in{\rm Sym}$, which are related to $h_1,h_2,\dots$ by the
well--known formula
$$
\log H(t)=\sum_{n=1}^\infty \frac{p_n\,t^n}{n}\,.
$$
{F}rom this formula we get
$$
\psi_{\alpha,\beta}(p_1)=1, \qquad \psi_{\alpha,\beta}(p_n)=\sum
\alpha_i^n+(-1)^{n-1}\sum\beta_j^n \quad (n\ge2),
$$
and the continuity of $\psi_{\alpha,\beta}(p_n)$ is readily verified (see,
e.g., \cite{KOO}).
\endpf

Thus, every character $\psi_{(\alpha,\beta)}$ is, loosely speaking, a
M--mixture of countably many copies of $\psi_+$, countably many copies of
$\psi_-$ and a single copy of $\psi_{0,0}$ in proportions
$\alpha_1,\alpha_2,\ldots$,\, $\beta_1,\beta_2,\ldots$, and $\gamma$,
respectively. Here $\psi_{0,0}$ has all parameters equal to 0.

\par
This completes Kerov's construction. It actually provides all the characters,
although the proof of this completeness claim requires further techniques.

\begin{remark} Note that Kingman's boundary $\partial{\cal P}$ can be constructed in the
same way. In this case we consider $\rm Sym$ with  the basis of (augmented)
monomial symmetric functions. This basis is not preserved by the canonical
involution, and the operation of M--mixing involves only one elementary
character: the corresponding probability function is supported by the
partitions $(n)\in\cal P$, so that this character may be viewed as a
counterpart of $\psi_+$. As a consequence, the boundary is just $\Delta$. A
character with the probability function supported by $\{(1^n)\}$ still exists
but it plays the role of $\psi_{0,0}$, and not the role of $\psi_-$.
\end{remark}

\section{Finitary skeleton of the boundary $\partial\cal Z$}

Now we aim to adapt Kerov's construction  to the graph $\cal Z$ by dealing with
the algebra ${\rm QSym}$ and its basis $\{F_\lambda: \lambda\in\cal Z\}$ of
fundamental quasisymmetric functions. As before, we have two elementary
characters but the construction becomes more sophisticated because the
comultiplication is no longer cocommutative. Another complication is the lack
of a nice system of generators like $\{h_n\}$ or $\{p_n\}$. As a consequence,
we cannot write an analog of formula (\ref{psi}). On the other hand, the action
of the comultiplication on the $F$--functions has a relatively simple form
(compare (\ref{delta}) and (\ref{deltaS})), which enables us to control the
behaviour of characters on the elements of the basis rather than on a system of
generators.

\par Extending (\ref{delta}) we immediately obtain a simple expression for the
iterated comultiplication map $\delta^{(k)}:{\rm QSym}\to {\rm QSym}^{\otimes
k}$:
$$
\delta^{(k)}(F_\lambda)= \sum F_{\lambda(1)}\otimes\cdots\otimes F_{\lambda(k)}
$$
where the summation is over all splittings of zigzag diagram $\lambda$ at
arbitrary positions in $k$ consecutive zigzag subdiagrams
$\lambda(1)\sqcup\ldots\sqcup\lambda(k)$ some of which may be empty.

\par The definition of the endomorphisms $r_\omega$ and of the M--mixing
(Definition \ref{Mmix}) naturally extends to the algebra ${\rm QSym}$. However,
now the resulting character $\psi$ is sensitive to the ordering of
the couples $(\psi_i,\omega_i)$. Here is an explicit expression
for $\psi$:
\begin{equation}\label{mix}
\psi(F_\lambda)=\sum \,\omega_1^{\ell_1}\cdots
\omega_k^{\ell_k}\,\, \psi_1(F_{\lambda(1)})\cdots
\psi_k(F_{\lambda(k)})\,, \qquad\lambda\in {\cal Z}\,,
\end{equation}
where the summation is as above and $\ell_j=|\lambda(j)|$, so that
$\lambda(j)\in {\cal Z}_{\ell_j}$.

\par
We define the elementary characters $\psi_+,\,\psi_-: {\rm QSym}\to{\rm QSym}$
exactly as in Proposition \ref{elemchar}, the only change is that $\lambda$
ranges over ${\cal Z}$ and the $S$--functions are replaced by the
$F$--functions. Both proofs of Proposition \ref{elemchar} remain valid: in the
second proof we use the multiplication rule (\ref{mult}). Note that the
involution map (Definition \ref{inv}) interchanges $\psi_+$ and $\psi_-$.

\begin{proposition}
The elementary characters of $\,{\rm QSym}$ defined above are extensions of the
elementary characters of\/ ${\rm Sym}$.
\end{proposition}

\proof Consider the character $\psi_+$ of ${\rm QSym}$ and restrict it to ${\rm
Sym}$. Since $S_{(n)}=F_{(n)}$, we have $\psi_+(S_{(n)})=1$. If $\lambda\in\cal
Y$ is not a one--row diagram then the expansion of $S_\lambda$ in the
$F$--functions does not involve one--row zigzag diagrams, as follows from
\cite[Theorem 7.19.7, p. 361]{StanleyII}. Therefore, $\psi_+(S_\lambda)=0$.
Thus, the restriction of $\psi_+$ has the required form. For $\psi_-$, we use a
similar argument or apply the involution.
\endpf

\par
Now we can select arbitrary weights and apply the M--mixing to some number of
copies of $\psi_+$ and $\psi_-$ arranged in a sequence. It is convenient to
represent these data in the following form.

\begin{definition}\label{finpaint}
Define a {\it finitary oriented paintbox\/} to be a pair $(U_\up, U_\down)$ of
disjoint open subsets of the open interval $]0,1[$, each comprised of finitely
many interval components and such that their total Lebesgue  measure is ${\rm
Leb}(U_\up \cup U_\down)=1$. Let ${\cal U}^{(2)}_{\rm \,fin}$ denote the set of
all such pairs.
\end{definition}

Think of the interval components of $U_\up$ as having positive orientation, and
think of the interval components of $U_\down$ as having negative orientation.
With this convention, each interval component of  the union $U_\up \cup
U_\down$ is characterised by location,  length and orientation.
\par
We define M--mixture $\psi_{(U_\up,U_\down)}$ by assigning, to each interval
component $]a,b[\,\subset U_\up \cup U_\down\,$, the factor $\psi_+$ with
weight $b-a$ in the case $]a,b[\,\subset U_\up$, or the factor $\psi_-$ with
weight $b-a$ in the case $]a,b[\,\subset U_\down$; the order of factors is
maintained in accordance with the natural ordering of the intervals, that is, from
left to right.

\par Let
\begin{equation}\label{rank}
{\rm rank}:(U_\up,U_\down)\mapsto (\alpha,\beta)
\end{equation}
be the map assigning to each $(U_\up,U_\down)\in{\cal U}^{(2)}_{\rm \,fin}$ two
nonincreasing positive sequences $\alpha=(\alpha_1,\dots,\alpha_k)$ and
$\beta=(\beta_1,\dots\beta_l)$ of lengths of the interval components of
$U_\up$ and $U_\down$, respectively. Restricting $\psi_{(U_\up,U_\down)}$ to
${\rm Sym}$ we get the character $\psi_{\alpha,\beta}$ as defined in Section 4.
Conversely, given such $\alpha$ and
$\beta$ with $\sum\alpha_i+\sum\beta_j=1$, there
are as many finitary characters $\psi_{(U_\up,U_\down)}$ extending
$\psi_{\alpha,\beta}$ as there are distinct arrangements of these $k+l$
intervals within $]0,1[$, which is $(k+l)!$ in the case all lengths are
pairwise distinct.

\begin{proposition}
For a finitary oriented paintbox $(U_\up,U_\down)$ with $m$ interval
components, let $\omega_1,\ldots,\omega_m$ be the sequence of lengths of these
components  in the natural order, from left to right. Then for any
$\lambda\in {\cal Z}$
\begin{equation}
\label{UUF}
\psi_{(U_\up,U_\down)}(F_\lambda)=\sum \omega_1^{\ell_1}\cdots \omega_m^{\ell_m}
\end{equation}
where the sum expands over all splittings $\lambda=\lambda(1)\sqcup\dots\sqcup
\lambda(m)$ such that $\lambda(j)=(\ell_j)$ or $\lambda(j)=(1^{\ell_j})$, with
$\ell_j\geq 0$, according as $\omega_j$ is the length of   a $U_\up$--component
or a $U_\down$--component. In the case $\ell_j=0$ we understand $\lambda(j)$ as
the empty zigzag.
\end{proposition}

\proof This immediately follows from formula (\ref{mix}) and the definition of
$\psi_+$ and $\psi_-$.
\endpf

In the particular case $U_\down=\varnothing$ the formula for the probability
function corresponding to $\psi_{U_\up,U_\down}=\psi_{U_\up,\varnothing}$
appeared as \cite[Equation 6.4]{DFP}, this is just the specialisation of
$F_\lambda$ in the variables $\omega_j$ with no type distinction.

\begin{definition}[Oriented paintbox construction]\label{paintconstr}
We associate with finitary $(U_\up,U_\down)$  the following construction of a
random arrangement $\Pi=(\Pi_n)$. Informally, $\Pi$ breaks $\mathbb N$ into
blocks and puts them in some order, whereas the order within each of the blocks
coincides with either the standard order or with the inverse order. In more
detail, let $\xi_1,\xi_2,\ldots$ be an independent sample from the uniform
distribution on $[0,1]$. Keep in mind that the sample values are almost surely
distinct and that $\xi_j$ hits the complement to $U_\up\cup U_\down$ in $[0,1]$
with probability zero. Define a random permutation $\Pi_n$ of $[n]$ by placing
$i$ to the left of $j$ if one of the following conditions holds:
\begin{itemize}
\item $\xi_i$ and $\xi_j$ do not hit the same interval component of $U_\up\cup
U_\down$ and $\xi_i<\xi_j\,$;
 \item $\xi_i$ and $\xi_j$ hit the same interval component of $U_\up$ and $i<j\,$;
 \item $\xi_i$ and $\xi_j$ hit  the same interval component of $U_\down$ and $i>j\,$.
\end{itemize}
The number $n$ enters implicitly, in the form of the constraint
$n\geq\max(i,j)$, hence the algorithm  yields a coherent sequence of
permutations, {\it i.e.\/,} defines a random arrangement $\Pi=(\Pi_n)$ of
$\mathbb N$. Let us say that $\Pi$ is {\it directed\/} by $(U_\up,U_\down)$.
\end{definition}

In the case $U_\down=\varnothing$ the resulting permutation $\Pi_n$ is inverse
to a biased riffle shuffle introduced in \cite{DFP} and further studied in
\cite{Fulman}. In \cite[Definition 2.1]{Lalley}, a more general concept of
shuffle is introduced in terms of a measure--preserving mapping
$f:[0,1]\to[0,1]$. This  covers the general finitary $(U_\up,U_\down)$,
although that paper is really focussed on mappings $f$ increasing on its
intervals of continuity, which in our context still corresponds to
$U_\down=\varnothing$.
\par It is known
 that for Young's lattice there is also a way to transform an independent
sample into a Markov chain on $\cal Y$, but it is much more complicated as is
based on properties of the generalised RSK--correspondence \cite{VKSIAM}.

\vskip0.5cm

\begin{proposition} The distribution of the random arrangement
directed by $(U_\up,U_\down)$ is
$$
{\mathbb P}(\Pi_n=\pi_n)=\psi_{(U_\up,U_\down)}(F_\lambda)\,, ~~~{\rm
for~~}\lambda={\rm zs}(\pi_n),
$$
where $\psi_{(U_\up,U_\down)}(F_\lambda)$ is given by {\rm (\ref{UUF})}.
\end{proposition}

\proof Let $I_1,\ldots,I_m$ be the interval components  of the oriented
paintbox and let $\omega_1,\ldots,\omega_m$ be their lengths. For
$i=1,\dots,m$, let $\ell_i$ be the number of indices $j\in[n]$ such that
$\xi_j\in I_i$. If $I_i$ belongs to $U_\up$ then $I_i$ contributes an
increasing segment of length $\ell_i$ to $\Pi_n$, while if $I_i$ belongs to
$U_\down$ then the contribution is a decreasing segment of length $\ell_i$. By
convention, empty and singleton segments may be viewed both as increasing and
decreasing. Each splitting of $\pi_n$ in consecutive monotonic segments of
lengths $\ell_1,\ldots,\ell_m\geq 0$ contributes to ${\mathbb P}(\Pi_n=\pi_n)$
the probability
$$
{\mathbb P}(\xi_{\pi_n(1)},\ldots,\xi_{\pi_n(\ell_1)}  \in
I_1\,;~\ldots ~; \xi_{\pi_n(\ell_1+\ldots+\ell_{m-1}+1)}, \ldots
,\xi_{\pi_n(n)}\in I_m)
=\omega_1^{\ell_1}\cdots\omega_m^{\ell_m}
$$
provided the monotonicity type of each segment $\pi_n
(\ell_1+\ldots+\ell_{i-1}+1)$, \dots, $\pi_n(\ell_1+\ldots+\ell_i)$ agrees with
the orientation of $I_i$. Summing over all such splittings yields (\ref{UUF}).
\endpf

\begin{example}[$a$--shuffles \cite{ashuffle}] Formula (\ref{UUF}) simplifies
considerably in the case when $U_\downarrow=\varnothing$ and $U_\up$ is the
equispaced division of $[0,1]$ in $a$ intervals. In this case the permutation
is inverse to a riffle shuffle \cite{ashuffle} with $a$ piles and its
distribution is given by
$$
{\mathbb P}(\Pi_n=\pi_n)={n+a-k\choose n}a^{-n}\,,
$$
where $k-1$ is the number of descents in $\pi_n$. Reciprocally, with the roles
of $U_\up$ and $U_\down$ exchanged, the right--hand side gives the probability
${\mathbb P}(\Pi_n=\pi_n')$ for the reversed permutation. As $a\to\infty$ both
series converge to the distribution of a uniform permutation which has
$\prob(\Pi_n=\pi_n)=1/E_{n,k}$, where $E_{n,k}$ is the Eulerian number
\cite{StanleyI}. In \cite{Euler} we show that the distributions just described
are extreme within the smaller class of  coherent permutations which have the
{\it number} of descents as sufficient statistic. These distributions play for
$\cal Z$ a role comparable with that  of Ewens' distributions for $\cal P$, in
the sense that  Ewens' distributions are extreme among the distributions for
exchangeable partitions which have the number of blocks as sufficient statistic
\cite{Gibbs}. Still, there is a difference: Ewens' distributions can be
decomposed over $\partial{\cal P}$ (the mixing measure is known as the
Poisson--Dirichlet distribution \cite{CSP}).
\end{example}

\begin{example}
[`bi--interval' \cite{GN}]\label{bi} Suppose the oriented paintbox is composed
of one $\down$--interval $]0,\varphi[$ and one $\up$--interval $]\varphi,1[$
with common endpoint $\varphi\in \,]0,1[$. The random zigzag shape $X_n={\rm
zs}(\Pi_n)$ has the distribution supported by hook zigzags only, and its
probability function is
$$p(1^l, k+1)=\varphi^l(1-\varphi)^{k}\,,~~~k+l=n-1\,.$$
The distribution of the resulting arrangement $\lhd$ is also
easy to describe: with probability $\varphi$ the integer $n$ is
$\lhd$--ordered below $[n-1]$ and with probability $1-\varphi$
it is $\lhd$--ordered above $[n-1]$. Thus the initial ranks
$\rho_n:=\Pi_n^{-1}(n)$ are independent, assume only extreme
values and have a representation
$$\rho_n=1\cdot {\bf 1}(\eta_n=0)+n\cdot{\bf 1}(\eta_n=1),$$
where $\eta_n$'s are independent Bernoulli random variables with
${\mathbb P}(\eta_n=0)=\varphi$.
\end{example}

\section{The general oriented paintbox}

\subsection{The closure}
The set ${\cal U}^{(2)}_{\rm \,fin}$ introduced in Definition \ref{finpaint} is
an analog of $\Delta^{(2)}_{\rm fin}$.  We seek now for an analog of
$\Delta^{(2)}$ --- a  topological space which would serve as a set of
parameters for the closure of the set of finitary characters. With reference to
\cite{RCS} the logical choice of that space is
$$
{\cal U}^{(2)} =\{(U_\up,U_\down): \text{$U_\up, U_\down$ are open disjoint
subsets of $]0,1[$}\}.
$$
Let us emphasise that, by definition, $U_\up$ and $U_\down$ are contained in
the open interval $]0,1[$ hence include neither $0$ nor $1$. We endow ${\cal
U}^{(2)}$ with a kind of Hausdorff distance, defined for
$(U_\up,U_\down),~(V_\up,V_\down)\in {\cal U}^{(2)}$ to be the infimum $\theta$
such that the $\theta$--inflation of the complementary closed set $U_\up^{\rm
compl}:=[0,1]\setminus U_\up$ covers $V_\up^{\rm compl}$, the
$\theta$--inflation of $V_\down^{\rm compl}$ covers $U_\down^{\rm compl}$, and
the same holds with the roles of $(U_\up,U_\down)$ and $(V_\up,V_\down)$
swapped.

\par The convergence of a sequence $((V_\up(j),V_\down(j)), j=1,2,\ldots)$ to some
$(U_\up,U_\down)$ in this metric is described by the following rule. Introduce
a notation: for any $\epsilon>0$, let $(U^\epsilon_\up,U^\epsilon_\down)$ stand
for the paintbox obtained from $(U_\up,U_\down)$ by deleting all interval
components with length less or equal to $\epsilon$. Then, for arbitrarily small
fixed $\epsilon$, different of all interval lengths represented in $U_\up\cup
U_\down$, for $j$ large enough, $V_\up^\epsilon(j)\cup V_\down^\epsilon(j)$
must have the same number of intervals as $U_\up^\epsilon\cup
U_\down^\epsilon$, with the same orientation pattern. Moreover, the record of
the endpoints of all intervals for $(V_\up^\epsilon(j),V_\down^\epsilon(j))$
must approach  the similar record for $(U_\up,U_\down)$, as $j\to \infty$.

\par For instance, $(V_\up(j),V_\down(j))\to (\varnothing, \varnothing)$ in ${\cal
U}^{(2)}$ if and only if the length of the longest interval component
approaches $0$.

\begin{proposition}\label{dense}
  ${\cal U}^{(2)}_{\rm \,fin}$ is dense in ${\cal
U}^{(2)}$.
\end{proposition}

\proof Let $(U_\up,U_\down)\in{\cal U}^{(2)}$ be arbitrary. An approximating
sequence
$$
\{(V_\up(j),V_\down(j))\}\to(U_\up,U_\down)
$$
of finitary paintboxes can be obtained as follows. Choose a sequence
$\epsilon_j\downarrow 0$. Given $j$, first take $(U_\up^{\epsilon_j},
U_\down^{\epsilon_j})$ and then extend it to a finitary paintbox by breaking
the complement $[0,1]\setminus (U_\up^{\epsilon_j}\cup U_\down^{\epsilon_j})$
in finitely many intervals of length not exceeding $\epsilon_j$.
\endpf

\begin{proposition}\label{compact}
${\cal U}^{(2)}$ is a compact space.
\end{proposition}

\proof Let us show that any sequence $((V_\up(j),V_\down(j)), j=1,2,\ldots)$ contains a
converging subsequence. For any fixed $\epsilon>0$, the total number of
interval components in $(V_\up^\epsilon(j),V_\down^\epsilon(j))$ does not
exceed $1/\epsilon$. Therefore, passing to a subsequence of indices $j$, we may
assume that the total number of intervals in
$(V_\up^\epsilon(j),V_\down^\epsilon(j))$ together with their orientation
pattern stabilise as $j\to\infty$, and moreover, the record of the endpoints of
the intervals is converging. Furthermore, fix $\epsilon_1>\epsilon_2>\dots$
converging to 0. Then, using Cantor's diagonal process we can choose a further
subsequence of indices $j$ such that the above convergence property holds for any
$\epsilon_j$. It is readily seen
that the resulting subsequence of paintboxes
is indeed converging to some element in  ${\cal U}^{(2)}$.
\endpf

\par
The oriented paintbox construction introduced before for finitary paintboxes
(Definition \ref{paintconstr}) extends literally to arbitrary
 $(U_\up,U_\down)\in {\cal U}^{(2)}$.
We were careful to formulate the first rule to make it working in the more
general situation: when $\xi_i,~\xi_j$ do not hit the same component interval,
one of these points may fall in the complement to $U_\up\cup U_\down$ in which
case the order is maintained according as $\xi_i<\xi_j$ or $\xi_i>\xi_j$. A
somewhat different description of  the construction is given in \cite{JW}.

\begin{proposition}\label{contin}
The correspondence $(U_\up,U_\down)\mapsto \Pi=(\Pi_n)$ between general
oriented paintboxes and random arrangements provided by the extension of
{\rm \,Definition \ref{paintconstr}} is continuous.
\end{proposition}

\proof The assertion  means that for any $n=1,2,\dots$ and any permutation
$\pi_n$ of $[n]$, the probability ${\mathbb P}(\Pi_n=\pi_n)$ is a continuous
function of $(U_\up,U_\down)$. Recall that the random permutation $\Pi_n$ is
uniquely determined by $(U_\up,U_\down)$ and the first $n$ sample values
$\xi_1,\ldots,\xi_n$. Observe that, for small $\epsilon>0$, the probability
that $|\xi_i-\xi_j|\le\epsilon$ for at least one pair  of distinct
indices $i,j\in [n]$ is $O(\epsilon)$. It follows that the event
$\Pi_n=\pi_n$ is determined, up to an event of probability $O(\epsilon)$,
by the mutual location of the sample
$\xi_1,\dots,\xi_n$ and $(U_\up^\epsilon, U_\down^\epsilon)$ only.
This easily implies the continuity assertion.
\endpf

\begin{corollary}\label{correspondence}
Let\/ $\Pi=(\Pi_n)$ be the random arrangement corresponding to an arbitrary
$(U_\up,U_\down)\in{\cal U}^{(2)}$. Then each probability of the form ${\mathbb
P}(\Pi_n=\pi_n)$ depends on the zigzag shape ${\rm zs}(\pi_n)$ only, so that
$\Pi$ determines a probability function $p=p_{(U_\up,U_\down)}\in\PF$.
Moreover, the corresponding functional $\psi_{(U_\up,U_\down)}$, given by
$\psi_{(U_\up,U_\down)}(F_\lambda)=p_{(U_\up,U_\down)}(\lambda)$, is
multiplicative and hence is a character. The set of characters obtained in this
way is the closure of the set of characters derived from the elements of ${\cal
U}^{(2)}_{\,\rm fin}$.
\end{corollary}

\proof This immediately follows from Propositions \ref{dense} and \ref{contin},
and for the last claim of the corollary, we also use Proposition \ref{compact}.
\endpf

Observe that the projection (\ref{rank}) extends, in a natural way, to a
continuous projection
$$
{\rm rank}: {\cal U}^{(2)}\to\Delta^{(2)}, \qquad
(U_\up,U_\down)\mapsto(\alpha,\beta),
$$
where the coordinates of $\alpha$ and $\beta$ are the lengths of the interval
components in $U_\up$ and $U_\down$, respectively, written in decreasing order.
Recall also that there is a natural projection $\partial\cal Z\to\partial\cal
Y$, see Proposition \ref{projection}.

\begin{corollary} Under the projection $\partial{\cal Z}\to \partial{\cal Y}$,
a character $\psi_{(U_\up,U_\down)}$ with $(U_\up,U_\down)\in{\cal U}^{(2)}$ is
sent to the character $\psi_{\alpha, \beta}$ of\/ ${\rm Sym}$ with parameters
$(\alpha,\beta)={\rm rank}(U_\up,U_\down)$.
\end{corollary}

\proof Since the embedding ${\rm Sym}\hookrightarrow{\rm QSym}$ respects the
comultiplication, the claim is true for finitary characters. The general case
follows by continuity of the projection.
\endpf

\begin{definition}\label{iota}
For $n=2,3,\dots$ define an embedding $\iota_n:{\cal Z}_n \hookrightarrow{\cal
U}^{(2)}\subset\partial{\cal Z}$ by assigning to each $\lambda\in{\cal Z}_n$ a
finitary paintbox $\iota_n(\lambda)=(V_\up,V_\down)$ which mimics the zigzag,
as follows. Recall that zigzags $\lambda\in\cal Z$ can be encoded by words
$w(\lambda)$ in the binary alphabet $\{+,-\}$, see the discussion following
Definition \ref{zigzag}. A word $w(\lambda)$ is further transformed into a
paintbox $(V_\up,V_\down)\in {\cal U}^{(2)}_{\,\rm fin}$, with the interval
components of $V_\up$ corresponding to the plus--clusters of $w(\lambda)$, and
those of $V_\down$ corresponding to the minus--clusters, taking the interval
lengths to be proportional to the cluster lengths, in the same order. For
example, $w(4, 1^2, 3)=+++---++\,$, and its image $\iota_{\,9}\,(4,1^2,3)
=(V_\up,V_\down)$ has $V_\up=\,]0,{\frac38}[\,\cup\,]{\frac68},{\frac88}[$ and
$V_\down=\,]{\frac38},{\frac68}[\,$.
\end{definition}

The next result is a law of large numbers for random arrangements derived from
oriented paintboxes.

\begin{proposition}\label{LLN}
Fix an oriented paintbox $(U_\up, U_\down)\in{\cal U}^{(2)}$ and
let $\Pi=(\Pi_n)$ be the corresponding random arrangement as
defined above. Write $\Pi$ as a random standard path
$(\lambda_0\nearrow\lambda_1\nearrow\dots)$ in ${\cal Z}$. Then
almost surely $\iota_n(\lambda_n)$ converges to $(U_\up,
U_\down)$ in the ambient space ${\cal U}^{(2)}$.
\end{proposition}

\proof Let $(V_\up(n),V_\down(n))=\iota_n(\lambda_n)$. We have to prove that
almost surely $(V_\up(n),V_\down(n))\to(U_\up,U_\down)$ as $n\to\infty$. We
will use the description of the convergence given in the beginning of this
section. Fix an arbitrarily small $\epsilon>0$ different from the lengths of
intervals in the paintbox $(U_\up,U_\down)$. Let $]a,b[$ be an interval
component in $U_\up$ with $b-a>\epsilon$; we assert that there is a similar
interval in $V_\up(n)$ for all $n$ large enough. Indeed, let
$(\xi_1,\xi_2,\dots)$ be the i.i.d random variables entering the paintbox
construction. Fix an arbitrary index $j$ such that $\xi_j$ hits $]a,b[\,\subset
U_\up$ (it exists almost surely). By the strong law of large numbers applied to
$(\xi_m)$, the proportion of natural numbers $m\in[n]$ placed in $\Pi_n$ to the
left of $j$ converges to $a$, and the length of the increasing run containing
$j$ is asymptotic to $(b-a)n$. It follows that the position (on the
$1/(n-1)$--scale) of the plus--cluster in $w(\lambda_n)$ corresponding to this
run is close to the interval $]a,b[$.

\par Conversely, each plus--cluster in $w(\lambda_n)$,
longer than $\epsilon n$, must correspond to some
$\up$--interval. This is argued by noting that the probability
that an increasing run longer $\epsilon n$  stems from more than
one such interval approaches zero as $n\to\infty$, in
consequence of the fact that for any interval the maximal index
$j\leq n$ with $\xi_j$ hitting the interval is likely to satisfy
$n-j=o(n)$, while the minimal $j$ with $\xi_j$ hitting the
interval is just a finite random variable. A similar argument
works for $\down$--intervals and minus--clusters.
\endpf

\begin{corollary}\label{injectivity}
The correspondence $(U_\up,U_\down)\mapsto p_{(U_\up,U_\down)}$ defined in
{\rm \,Corollary \ref{correspondence}} is injective.
\end{corollary}

Let us  summarise the main results of this section:

\begin{proposition}\label{summary}
The correspondence $(U_\up,U_\down)\mapsto p_{(U_\up,U_\down)}$
afforded by the oriented paintbox construction is a
homeomorphism of the compact space ${\cal U}^{(2)}$ onto a part
of the boundary $\partial\cal Z$.
\end{proposition}

\proof Indeed, this follows from Corollaries \ref{correspondence} and
\ref{injectivity}.
\endpf

In Section 7 we will show that the image of ${\cal U}^{(2)}$ actually coincides
with the whole boundary $\partial\cal Z$.

\subsection{Mixing}

Here mixing  means exploiting the oriented paintbox construction with a {\it
random\/} $(U_\up,U_\down)$ chosen from some probability distribution on ${\cal
U}^{(2)} $. Proposition \ref{summary} obviously extends to the mixed case,
hence yields a continuous injection of the space of probability distributions
on ${\cal U}^{(2)}$ with weak topology into the space $\PF$ of probability
functions on $\cal Z$.

\begin{example}[a random bi--interval]
Suppose  as in Example \ref{bi} that $]0,1[$ is divided  in one
$\down$-- and one $\up$--interval as $]0,\varphi[\,\cup\, ]\varphi,1[$, but
this time $\varphi$ is random with a beta density ${\mathbb P}(\varphi\in {\rm
d}x)= x^{\theta_1-1}(1-x)^{\theta_2-1}{\rm d}x/{\rm B}(\theta_1,\theta_2)$,
where $\theta_1,\theta_2>0$. The sequence of initial ranks $(\rho_n)$ is no
longer independent. Rather, the number of $1$'s in the sequence is a sufficient
statistic: the conditional probability of $\rho_n=1$ given
$\rho_1,\ldots,\rho_{n-1}$ is $(k+\theta_1)/(k+\ell+\theta_1+\theta_2)$ where
$k$ is the number of $1$'s among $\rho_1,\ldots,\rho_{n-1}$ and $\ell=n-1-k$.
An astute reader might have noticed that in the representation
$$\rho_n=1\cdot {\bf 1}(\eta_n=0)+n\cdot{\bf 1}(\eta_n=1)$$
the $\eta_n$'s are exchangeable Bernoulli variables that follow
Polya's urn model \cite{Feller}.
\end{example}

\subsection{Quasi--uniform distributions}
For an interval component $]\varphi,b[$ of $U_\up$ we think
of the left endpoint $\varphi$ as the initial point, and of $b$
as the terminal point of the interval.
 For an interval $]a,\varphi[$ of $U_\down$ the convention about the endpoints is
 reversed.
With each $(U_\up,U_\down)\in {\cal U}^{(2)}$ we associate a
probability measure $\nu$ on $[0,1]$ obtained by weeping out the
Lebesgue measure of each interval component of $U_\up\sqcup
U_\down$ to its initial point:
$$
\nu=\sum_{]\varphi,b[\,\subset U_\up}\delta_\varphi +
\sum_{]a,\varphi[\,\subset U_\down}\delta_\varphi+{\rm Leb}|_{[0,1]\setminus
(U_\up \cup U_\down)}.
$$
Such measures were dubbed {\it quasi--uniform\/} in \cite{JW}. In the case
$U_\up=\varnothing$ quasi--uniform measures coincide with {\it uniformised}
measures introduced in \cite{RCS}. By \cite[Lemma 3.3]{JW}, the quasi--uniform
measures can be characterised by the property that the inequalities
$\nu[0,x[\,\leq x\,\leq \nu[0,x]$ hold for $\nu$--almost all $x$.

\begin{proposition}
The above correspondence $(U_\up,U_\down)\mapsto\nu$ is a homeomorphism between
the space ${\cal U}^{(2)}$ with the Hausdorff--type metric as defined in
{\rm \,Section 6.1}, and the space of quasi--uniform measures endowed with the weak
topology.
\end{proposition}

\proof It is a straightforward check that the metric agrees with the Skorohod
topology on the space of c{\`a}dl{\`a}g distribution functions for
quasi--uniform measures.
\endpf

\par
The role of this class of distributions stems from the following observation.
Suppose $\nu$ is a quasi--uniform distribution corresponding to some
$(U_\up,U_\down)$. Let $\xi_j$'s be independent uniformly distributed on
$[0,1]$ and $\varphi_j$ be the initial point of the $(U_\up,U_\down)$--interval
discovered by $\xi_j$, with the convention $\varphi_j=\xi_j$ in case $\xi_j$
falls in $[0,1]\setminus (U_\up\cup U_\down)$. Then the probability law of
$\varphi_j$ is $\nu$. The value $\varphi_j$ may be defined intrinsically in
terms of the induced order $\lhd$ as the frequency of integers ordered below
$j$,
$$\varphi_j=\lim_{n\to\infty}\# \{1\leq i\leq n:~ i\lhd j\}/n\,,$$
which suggests to call $\varphi_j$ the {\it height} of $\lhd$ at $j$.

\subsection{An alternative description}

The sequence of heights alone determines only a {\it weak order}
on $\Nat$ by the formula \eq\label{CS} i\precsim
j\Longleftrightarrow \varphi_i\leq \varphi_j\,. \en The relation
$\precsim$ projects to a total order on $\Nat/\sim$ when
factored by the equivalence relation \eq\label{PS} i\sim j
\Longleftrightarrow \varphi_i= \varphi_j\,. \en
 The random partition of $\Nat$ into some collection of blocks
associated with $\sim$ is {\it exchangeable,\/} {\it i.e.\/} its
probability law is  invariant under all bijections
$\Nat\to\Nat$, and the realisation of $\sim$ via (\ref{PS}) is
an instance of Kingman's  paintbox representation of extreme
exchangeable equivalence relations on $\Nat$, also called {\it
partition structures} \cite{Ki1, Ki2}. Similarly, the
realisation of $\precsim$ by (\ref{CS})  is a version of the
representation of extreme exchangeable ordered partitions,  also
called {\it composition structures\/} \cite{RCS, GP}.

\par
The  strict order $\lhd$ extends the weak order $\precsim$ (off
the diagonal $\{(i,i):~i\in \Nat\}$)  by arranging in a certain
way the integers within each block of $\sim$. The modus operandi
of this extension is the following. Divide the collection of
intervals making up $U_\up\cup U_\down$ into {\it bi--intervals}
$]a,\varphi[\cup]\varphi,b[$ for each initial point $\varphi$.
If $\varphi$ is the initial point of a single interval, either
$]a,\varphi[\subset U_\down$ or $]\varphi,b[\subset U_\up$, we
adopt the convention that another component is empty, with
$b=\varphi$ or $a=\varphi$, respectively. The correspondence
between bi--intervals and initial points is bijective. It is
seen from the oriented paintbox construction that for a generic
initial point $\varphi$ the equality $\varphi_n=\varphi$ means
that $\xi_n$ hits the bi--interval $]a,\varphi[\cup]\varphi,b[$
corresponding to $\varphi$, thus those $n$'s comprise a
nonsingleton block of $\sim$, call it $B_\varphi$. By the
construction, in the case $\xi_n$ hits $]a,\varphi[$ the integer
$n$ is $\lhd$--ordered below all smaller $j$'s with $j\in
B_\varphi$, while in the case $\xi_n\in ]\varphi,b[$ the integer
$n$ is $\lhd$--ordered above all smaller $j$'s with $j\in
B_\varphi$. Write $s_n=\,\downarrow$ in the first case,
$s_n=\,\uparrow$ in the second case, and write $s_n=\cdot$ if
$\xi_n$ misses $U_\up\cup U_\down$. All this encodes $\lhd$ into
a random sequence $((\varphi_n,s_n), n=1,2,\ldots)$ by the rule:
\eq\label{O} ~~(i<j)\&(i\lhd j)~~\Longleftrightarrow
~~\varphi_i<\varphi_j~~{\rm or}~~ (\varphi_i=\varphi_j)\&
(s_j=\,\up)\,. \en For $n=1,2,\ldots$ the pairs
$(\varphi_j,s_j)$ are independent and identically distributed,
with $\varphi_n$ having a quasi--uniform distribution $\nu$ and
the conditional law of $s_n$ given by
$$
{\mathbb P}(s_n=\,\uparrow\giv \varphi_n=\varphi)
=\frac{b-\varphi}{b-a}\,, \qquad {\mathbb P}(s_n=\,\down\giv
\varphi_n=\varphi) =\frac{\varphi-a}{b-a}
$$
for each $\varphi\in[0,1]$ being the initial point  of a bi--interval
$]a,\varphi[\,\cup\,]\varphi,b[$, and ${\mathbb P}(s_n=\,{\bf\cdot}\,\giv
\varphi_n=\varphi)=1$ for each $\varphi\in [0,1]\setminus(U_\up\cup U_\down)$.

\par
The representation (\ref{O})  of $\lhd$ involving the heights is
canonical; other representations may be obtained by replacing
$\varphi_n$'s by $\theta_n=f(\varphi_n)$ with
$f:[0,1]\to{\mathbb R}$ strictly increasing on the support of
$\nu$. From a statistical perspective, extending the relation
$i\precsim j\Longleftrightarrow \theta_i\leq\theta_j$ may be
regarded   as `breaking  ties' in the sequence with repetitions
$(\theta_n)$, and this extension is very natural since it agrees
with spreadability  of the resulting total order  (the
concept of speadability is discussed in the next section).

\vskip0.5cm

\begin{remark}
For $i,j\in B_\varphi$ the number of integers $\lhd$--ordered
between $i$ and $j$ has a binomial distribution. Moreover, no integer
$n>\max(i,j)$ can be $\lhd$--ordered between $i$ and $j$. Each $B_\varphi$ is
therefore a saturated $\lhd$--chain with the property that there are only
finitely many elements between any $i,j\in B_\varphi$.
\end{remark}

\section{Completeness of the description of the boundary}

The random order $\lhd$ induced by an oriented paintbox is not exchangeable,
since its distribution may be affected by bijections $\Nat\to\Nat$ (with the
sole exception $U_\up\cup U_\down= \varnothing$ when each $\Pi_n$ is a uniform
random permutation). Still, the order $\lhd$ possesses a weaker symmetry
property of {\it spreadability\/}.
 Spreadability is the key
property which allows reducing description of all random arrangements $\Pi$
with descent sets as sufficient statistics to a recent result by Jacka and
Warren \cite{JW}. See \cite{Aldous} for a survey of earlier results on
spreadability and its relation to exchangeability, and see \cite{Kallenberg}
for an updated exposition (in \cite{Kallenberg} `spreadability' is also called
`contractability').

\par
For a subset $B\subset \Nat$ let us call {\it ranking}  the standard order
isomorphism of $B$ onto some $[n]$ or onto $\Nat$, thus the ranking sends $b\in
B$ to $j$ when $b$ is the $j$th smallest element in $B$. For $\pi_n$ a
permutation of $[n]$ and $j\in[n]$, let $\tau_j(\pi_n)$ be the permutation of
$[n-1]$ obtained by removing $j$ from the sequence $\pi_n(1),\cdots,\pi_n(n)$
and replacing the elements in the reduced sequence by their ranks.

\begin{proposition}\label{L}
Suppose $\Pi_n$ is a random permutation of\/ $[n]$ which satisfies
\begin{equation}\label{suff}
{\mathbb P}(\Pi_n=\pi_n)=p({\rm zs}(\pi_n))
\end{equation}
for some function $p:{\cal Z}_n\to [0,1]$. Then  the distribution of
$\tau_j(\Pi_n)$ is the same for all $j=1,\dots,n$ and satisfies the counterpart of
{\rm (\ref{suff})} with the function $p$ extended to ${\cal Z}_{n-1}$ by the
virtue of recursion {\rm (\ref{recY})} for the graph $\cal Z$.
\end{proposition}

\proof It is sufficient to consider the case when $\Pi_n$ is uniformly
distributed on the set of permutations with some fixed zigzag shape $\lambda$,
because the general case is a mixture.  Given $j=2,\dots,n$, divide the set of
such permutations in two disjoint subsets $A$ and $B$, with $A$ being the set
of permutations in which $j-1$ and $j$ occupy adjacent positions, and $B$ being
all other permutations which have $j-1$ and $j$ separated.
 We have $\tau_{j-1}(A)=\tau_j(A)$, meaning the equality of multisets, just
 because
$\tau_{j-1}(\pi_n)=\tau_j(\pi_n)$ for $\pi_n\in A$. The latter does not hold
for $B$, but exchanging $j$ and $j-1$ preserves the zigzag shape and, moreover,
yields a bijection of $B$, thus also $\tau_{j-1}(B)=\tau_j(B)$. Therefore
$\tau_j(\Pi_n)\ed \tau_{j-1}(\Pi_n)$ (equality of distributions), hence by
induction all distributions $\tau_j(\Pi_n)$ are the same for $1\leq j\leq n$.
The second assertion is true for $j=n$ (by the  definition of graph $\cal Z$
and (\ref{suff})) hence, from the above, it is true for all $j\leq  n$.
\endpf

\par
For $\Pi$ a random arrangement of $\Nat$ and an infinite set $B\subset\Nat$ let
$\Pi|_{B}$ denote the random arrangement obtained by ranking $B$, and for a
finite set $B\subset \Nat$ with $\#B=n$ let $\Pi|_{B}$ be the random
permutation of $[n]$ obtained in the same way. Call $\Pi$ {\it spreadable\/} if
$\Pi|_B$ has the same probability distribution as $\Pi$ for every infinite
$B\subset\Nat$. For $\Pi=(\Pi_n)$ the identity $\Pi|_{[n]}=\Pi_n$ holds by
definition, and in view of the fact that the law of $\Pi$ is completely
determined by the laws of coherent permutations $\Pi_n$, the spreadability
property just says that $\Pi|_{B}\ed\Pi_n$ for each $B\subset\Nat$ with $\#B=n$
and each $n$. This definition of spreadability for $\lhd$ agrees with that for
the random array $Z_{i,j}={\bf 1}(i\lhd j)$ (see \cite{Aldous}), meaning
$(Z_{i,j})\ed (Z_{T(i),T(j)})$ for all increasing $T:\Nat\to\Nat$. The next
corollary follows by iterated application of Proposition \ref{L}. \vskip0.5cm

\begin{corollary}\label{reduce}
If random permutations $(\Pi_n)$ are  coherent and each\/ $\Pi_n$ is
conditionally uniform given ${\rm zs}(\Pi_n)$, then the arrangement
$\Pi=(\Pi_n)$ is spreadable.
\end{corollary}
\noindent

The above definition of spreadability works equally well for the
set $\mathbb Z$. Note also that the oriented paintbox
construction can be trivially modified to produce random
arrangements of $\mathbb Z$, we just need to label
 the random variables $\xi_j$ by $\mathbb Z$
instead of $\Nat$.

\begin{proposition}\label{JW}
{\rm (Jacka and Warren \cite[Theorem 3.4]{JW})} Each spreadable
random arrangement on $\mathbb Z$ has the same probability law
as the one defined by the oriented paintbox $(U_\up,U_\down)$
chosen at random from some probability distribution on ${\cal
U}^{(2)}$.
\end{proposition}

\par
To apply Proposition \ref{JW} it remains to show that

\begin{lemma}\label{extend} Each spreadable arrangement on $\Nat$ has a
distributionally unique extension to a spreadable arrangement on
$\mathbb Z$.
\end{lemma}

\proof  For $\lhd$ a spreadable arrangement on $\mathbb Z$ the distribution is
uniquely determined by the distributions of $\lhd|_{[-n, n]}$ for
$n=1,2,\ldots$. On the other hand, for $\lhd$ a spreadable arrangement on
$\Nat$ the restriction $\lhd|_{[2n+1]}$ can be uniquely transported to
$[-n,n]$. It is easy to check that the pushforwards are consistent for various
values of $n$, hence determine an arrangement of $\mathbb Z$.
\endpf

\noindent
Putting the things together we arrive at our principal  conclusion.

\begin{theorem}\label{main}
The boundary $\partial{\cal Z}$ of the graph of zigzag diagrams is homeomorphic to
${\cal U}^{(2)}$. Each coherent sequence of permutations $\Pi$
with descent set as a sufficient statistic can be represented by
the oriented paintbox construction with some random
$(U_\up,U_\down)\in {\cal U}^{(2)}$. The distribution of the
paintbox  $(U_\up,U_\down)$ representing a given  $\Pi$ is
unique.
\end{theorem}

\proof The second claim follows from Corollary \ref{reduce}, Lemma
\ref{extend}, and Proposition \ref{JW}. Then the first claim becomes evident
due to Proposition \ref{summary}. The third claim follows by either recalling
that the set $\PF$ of probability functions is a Choquet simplex or, more
constructively, by appealing to Proposition \ref{LLN} which allows to identify
a random oriented paintbox with the almost sure limit of zigzags ${\rm
zs}(\Pi_n)$, suitably embedded in ${\cal U}^{(2)}$.
\endpf

\section{Complements and open questions}

\subsection{The involution}
The conjugation of zigzags
$\lambda\mapsto\lambda'$, defined
and interpreted in terms of binary encoding
in Section 2,
differs from another natural involutive operation on zigzags, see
\cite[Exercise 7.94]{StanleyII}. In terms of descent sets $D\subset[n-1]$, the
latter amounts to the operation $D\mapsto[n-1]\setminus D$, which is the same
as just switching $+\leftrightarrow-$. In terms of ribbon Young diagrams
$\bar\lambda$ (see Section 2), our version is related to the conventional
symmetry of (skew) Young diagrams, while the alternative version is related to
the reflection about the bisectrix of the first quadrant. An advantage of our
definition is that our operation induces a natural involution
$\pi_n\mapsto\pi'_n$ on the standard paths in the graph $\cal Z$ (in terms of
infinite paths this means inversion of total order on $\Nat$), while another
version leads to a cumbersome transformation of paths.
\par
These two operations on zigzags determine two distinct
involutive automorphisms of $\rm QSym$, $\rm inv$ and
$\widehat\omega$, which, however, coincide on $\rm Sym$ with the
canonical involution $\omega$ of that algebra. For
$\widehat\omega$, a simple proof of this fact is given in
\cite[solution to Exercise 7.94]{StanleyII}, and for $\rm inv$
the argument is similar.

\subsection{The empty paintbox}

The paintbox $U_\up=U_\down=\varnothing$ corresponds to the only
{\it exchangeable} random arrangement
$\Pi=(\Pi_n)$ with each $\Pi_n$ being the uniform random permutation of $[n]$.
Thus,
$$
p_{(\varnothing,\varnothing)} (\lambda_n)=d(\lambda_n)/n! \,\,~~~~\lambda_n\in
{\cal Z}_n\,.
$$
The corresponding character $\psi_{(\varnothing,\varnothing)}$ projects to the
Plancherel character $\psi_{0,0}$ of ${\rm Sym}$. Thus,
$\psi_{(\varnothing,\varnothing)}$ may be viewed as an analog of the Plancherel
character. There is a large literature on properties of the descent set for
uniform permutations, {\it e.g.\/}, it is well known that in this  case the
distribution of the number of descents is close to normal with mean $(n-1)/2$,
see \cite{Bender}. It would be interesting to learn which results on the
Plancherel measure on the Young diagrams and the corresponding measure on the
paths in the Young graph $\cal Y$ translate to the graph $\cal Z$, see
\cite{Oshanin} for some work in this direction.

\subsection{The Martin boundary}
The concept of Martin boundary for graded graphs \cite{KOO} may
be well adapted to our context; this corresponds to the entrance
boundary \cite{Kemeny} for the inverse Markov chain. Given two
zigzags $\mu,\lambda\in {\cal Z}$ let $d(\mu,\lambda)$ denote
the number of paths (if any) $\mu\nearrow\dots\nearrow\lambda$
in ${\cal Z}$ ascending from $\mu$ to $\lambda$ and set
\begin{equation}\label{martin}
K(\mu,\lambda)=\frac{d(\mu,\lambda)}{d(\lambda)}\,.
\end{equation}
For $\lambda_n\in {\cal Z}_n$ fixed, the function $\mu\mapsto K(\mu,\lambda_n)$
satisfies the recursion (\ref{recY}) in $\mu$ up to level $n-1$. A sequence
$(\lambda_n\in{\cal Z}_n)$ is said to be {\it regular\,} if for any fixed
$\mu$ the numerical sequence $K(\mu,\lambda_n)$ has a limit. For a regular
sequence $(\lambda_n)$, the function
$$
p(\mu):=\lim_{n\to\infty} K(\mu,\lambda_n)
$$
is a solution to (\ref{recY}), hence an element of $\PF$. The {\it Martin
boundary\/} of ${\cal Z}$ is the collection of all $p\in \PF$ which arise as
such limits. One can prove that the Martin boundary contains the boundary
$\partial{\cal Z}$, hence  $\partial{\cal Z}$ may be also called
the {\it minimal boundary\/}.
The inclusion relation between the two boundaries
also holds for the
general graded graphs, see for instance \cite[Chapter 1,
\S1]{KerovDiss}. Examples of graphs are known which have the Martin
boundary strictly larger than the minimal boundary
\cite{Gibbs}.

\begin{conjecture}\label{conj}
\begin{itemize}
\item[(i)] A sequence $(\lambda_n\in{\cal Z}_n)$ is regular if and only if the
points $\iota_n(\lambda_n)$
converge in the compactum ${\cal U}^{(2)}$.
\item[(ii)] $\iota_n(\lambda_n)\to (U_\up,U_\down)$ is equivalent to
$K(\mu,\lambda_n)\to \, p_{(U_\up,U_\down)}(\mu)$ for all $\mu\in \cal Z$.
\item[(iii)] The Martin boundary of the graph $\cal Z$ actually coincides with
its minimal boundary $\partial{\cal Z}={\cal U}^{(2)}$.
\end{itemize}
\end{conjecture}

The analogues of (i)--(iii) are known to be true for the Young graph $\cal Y$
\cite{KOO}, and for the graphs $\cal P$ and $\cal C$ \cite{Ki1, RCS,
KerovDiss}.

\subsection{The kernel $K(\mu,\lambda)$}
Consider the Young graph $\cal Y$ and the corresponding kernel
$K(\mu,\lambda)$ given by (\ref{martin}), where $\mu$ and
$\lambda$ are assumed to be Young diagrams. The algebraic
approach of \cite{KOO} to Thoma's theorem describing  the boundary
$\partial\cal Y$ relies on the following properties of the
kernel.

\begin{itemize}
\item The linear span (over $\mu$) of the set of functions
$$
f_\mu(\lambda)=n(n-1)\dots(n-m+1)K(\mu,\lambda), \qquad
\mu\in{\cal Y}_m, \quad \lambda\in{\cal Y}_n
$$
 is an algebra $A$ under pointwise multiplication.
\item This
algebra $A$ can be identified with the algebra $\rm Sym$ of
symmetric functions in such a way that
$$
f_\mu=S_\mu + {\rm lower~degree~terms}
$$
where $S_\mu$ is the Schur function.
\end{itemize}
\noindent Moreover, we have a very detailed information about the functions
$f_\mu$ viewed as elements of $\rm Sym$: these are certain `factorial'
analogues of Schur functions, called the {\it Frobenius--Schur functions\/},
see \cite{ORV}.

\par Similar approach applies to the graph $\cal P$: then the
corresponding functions $f_\mu$ are identified with natural factorial analogues
of monomial symmetric functions. As an application one gets an algebraic
derivation of Kingman's theorem about $\partial\cal P$, see \cite{KerovDiss}.
One can also prove similar claims for the graph $\cal C$: then $A$ is
identified with $\rm QSym$ and the functions $f_\mu$ turn into factorial
quasisymmetric monomial functions. Again, this yields an algebraic approach to
computing the boundary $\partial\cal C$, hence entails an alternative proof of
the main result of \cite{RCS}.

\begin{problem} Is it possible to extend this approach to the kernel
(\ref{martin}) on the graph $\cal Z$? \vskip0.5cm \noindent A positive answer
would lead to an alternative proof of Theorem \ref{main} and could also be
useful in studying Conjecture \ref{conj}.
\end{problem}

\subsection{Nonextreme solutions $p\in\PF$}

There has been an extensive study of concrete examples of {\it
nonextreme\/} solutions $p$ to recursion (\ref{recY})
for the graphs $\cal P$, $\cal C$, and $\cal Y$. These examples
arise from various probabilistic models and the problem of
harmonic analysis on the infinite symmetric group, see
\cite{RCS, CSP, BO}.

\begin{problem} Find natural examples of nonextreme solutions $p$ to
recursion (\ref{recY}) on the graph $\cal Z$.
\end{problem}

\subsection{The graph associated with the peak algebra}

A  graph $\cal Q$ closely related to ${\cal Z}$ can be viewed as a natural
`quasisymmetric analogue' of the so--called {\it Schur graph\/} ${\cal S}$
which appears in the theory of {\it projective\/} characters of the symmetric
groups and is associated with the classical $Q$--Schur functions \cite{Nazarov,
BorOls}. The standard finite paths in $\cal Q$ are in bijection with
permutations, like for $\cal Z$,
 but the set   of vertices is different. While a generic
vertex of $\cal Z$ corresponds to the descent set of a permutation, a vertex of
$\cal Q$ corresponds to the {\it peak set\/} of a permutation. Thus, the
analogue of recursion (\ref{recY}) for $\cal Q$ reflects the branching of
coherent permutations with the peak--set statistics. Alternatively, $\cal Q$
can be defined as the `multiplicative' graph associated with the {\it peak
algebra\/} $\rm Peak\subset\rm QSym$ and its distinguished basis introduced by
Stembridge \cite{Stembridge}.
\par Our method of constructing the boundary
can be extended to the graph $\cal Q$: here we use the fact that $\rm Peak$ is
closed under comultiplication. The algebra $\rm Peak$ has no natural
involution, hence in the construction of characters we may exploit only one
elementary character (as in the cases of ${\cal P}$ and ${\cal C}$). The
boundary of ${\cal Q}$ is parametrised by open sets $U\subset]0,1[\,$. The
completeness of the description of the boundary can be deduced from Theorem
\ref{main} using the fact that $\rm Peak$ can also be realised  as a quotient
of $\rm QSym$.
\par
The relations between the boundaries of the
graphs is represented by the diagram
$$
\begin{CD}
\partial{\cal Q}@>>>\partial{\cal Z} @<<< \partial{\cal C}\\
@VVV @VVV @VVV\\
\partial{\cal S} @>>>   \partial{\cal Y} @<<< \partial{\cal P}
\end{CD}
$$
Specifically, by the embedding $\partial{\cal P}
\hookrightarrow\partial{\cal Y}$ the parameter $\alpha\in
\Delta$ goes to $(\alpha,0)\in \Delta^{(2)}$. The boundary
$\partial{\cal S}$ is again $\Delta$ (as is read from
\cite{Nazarov}), but  the embedding  $\partial{\cal
S}\hookrightarrow\partial{\cal Y}$ sends $\alpha\in\Delta$ to
$(\alpha',\beta')\in \Delta^{(2)}$ with
$\alpha'=\beta'=(\alpha_1/2,\alpha_2/2,\ldots)$. Likewise,  the
embedding $\partial{\cal C}\hookrightarrow \partial{\cal Z}$
sends $U\in {\cal U}$ to $(U,\varnothing)\in {\cal U}^{(2)}$,
while by the embedding $\partial{\cal Q}\hookrightarrow
\partial{\cal Z}$ each $U\in {\cal U}$ is mapped to the pair
$(U_\up',U_\down')$ which is  obtained by halving each component
interval of $U$ and composing $U_\up'$ of right halves and
$U_\down'$ of left halves.

\subsection{The ring theorem}
Proposition \ref{P2} refers to the Kerov--Vershik `ring theorem'
\cite{VK-Roman, KerovDiss}. For reader's convenience we present
here a complete proof which follows the lines in
\cite{VK-Roman},
 but provides some details omitted in that paper.
Assume that
\begin{itemize}
\item[(i)] we are given a commutative unital algebra $A$ over
$\mathbb R$ together with a convex cone $K$ in $A$ (below we
write $a\ge b$ if $a-b\in K$),

\item[(ii)] the cone is generating ($K-K=A$) and stable under
multiplication ($K\cdot K\subseteq K$),

\item[(iii)] the unity $1$ in $A$ is also an order unity,
meaning $1\ge0$, and that for any $a\in K$ there exists a small
$\epsilon>0$ with $\epsilon a\le1$,

\item[(iv)] the cone $K$ is generated by a countable set of elements.
\end{itemize}
Denote by $P\subset A^*$ the set of all linear functionals
$\psi:A\to\mathbb R$ which are nonnegative on $K$ and normalised
by $\psi(1)=1$, and observe that $P$ is a convex set.

\begin{proposition}\label{ring}
Under these assumptions, a functional $\psi\in P$ is an extreme
point of $P$ if and only if $\psi$ is multiplicative, that is,
$\psi(ab)=\psi(a)\psi(b)$ for all $a,b\in A$.
\end{proposition}

\proof First of all, note that if $\psi$ vanishes on an element $a\ge0$ then
$\psi(ab)=0$ for any $b\in A$. Indeed, by (ii) and (iii), without loss of
generality, we may assume that $0\le b\le1$. Then we have
$$
0\le \psi(ab)\le\psi(ab)+\psi(a(1-b))=\psi(a)=0,
$$
whence $\psi(ab)=0$.
Assume now that $\psi$ is extreme. To  show that $\psi$ is
multiplicative, it suffices to prove that
$\psi(ab)=\psi(a)\psi(b)$ for $0\le a\le1$ and any $b$. Set
$a'=1-a$ and suppose first that both $\psi(a)$ and $\psi(a')$
are nonzero. Then the functionals
$$
\psi_a(b):=\frac{\psi(ab)}{\psi(a)} \quad {\rm and} \quad
\psi_{a'}(b):=\frac{\psi(a'b)}{\psi(a')}
$$
are well defined and are both elements of $P$. Since $\psi$ is a
convex combination of $\psi_a$ and $\psi_{a'}$ (with
coefficients $\psi(a)$ and $\psi(a')$, respectively), and $\psi$
is extreme, we have either $\psi=\psi_a$ or $\psi=\psi_{a'}$.
Consequently, we have either $\psi(b)=\psi_a(b)$ or
$\psi(b)=\psi_{a'}(b)$, and each of these relations implies
$\psi(ab)=\psi(a)\psi(b)$.
\par
Furthermore, if $\psi(a)=0$ then, as mentioned earlier,
$\psi(ab)=0$, which again implies
$\psi(ab)=\psi(a)\psi(b)$. Similarly, $\psi(a')=0$ implies
$\psi(a'b)=\psi(a')\psi(b)$, which is equivalent to
$\psi(ab)=\psi(a)\psi(b)$.
\par
Let us check the inverse implication. Denote by $P_{\rm ex}\subset P$ the
subset of extreme points. The convex set $P$ lies in the dual vector space
$A^*$ which we consider with the weak topology. By (iii) and (iv) there exists
a sequence $b_1, b_2,\ldots\in K$ spanning the cone and such that
$b_1\le1,b_2\le1,\dots$. By (ii), the family $\{b_1,b_2,\dots\}$ contains a
basis of $A$. It readily follows that the mapping
$\psi\mapsto(\psi(b_1),\psi(b_2),\dots)$ determines a homeomorphism of $P$ onto
a closed subset of the product space $[0,1]^\infty$, hence $P$ is a compact,
metrisable, and separable space. By Choquet's theorem, $P_{\rm ex}$ is a
$G_\delta$--subset in $P$, and each $\psi\in P$ is representable as a mixture
$$
\psi(a)=\int_{P_{\rm ex}}\varphi(a)\sigma(d\varphi)\qquad
\forall a\in A
$$
with some probability measure $\sigma$ on $P_{\rm ex}$.
\par
Using the embedding  $P\hookrightarrow [0,1]^\infty$  we may
view $\sigma$ as a probability measure on the product space
$[0,1]^\infty$, so that the coordinate functions in this space
become random variables. Now assume  $\psi$ multiplicative.
Then, in particular, $\psi(b_i^2)=(\psi(b_i))^2$ for every
$i=1,2,\dots$. It follows that each of the coordinate functions
has variance 0, hence is constant almost surely. This means that
$\sigma$ is supported by a single point, that is, $\psi\in
P_{\rm ex}$, as wanted.
\endpf

\vskip 0.5cm

\noindent  {\bf Acknowledgements.} We are grateful to Saul Jacka for making
accessible \cite{JW}.
Thanks go
to Natalia Tsilevich, Sergey Fomin and a referee for a number of helpful remarks.
The second author  was supported by the CRDF grant RUM1--2622--ST--04 and a
travel grant from the University of Utrecht.

\def\cprime{$'$} \def\polhk#1{\setbox0=\hbox{#1}{\ooalign{\hidewidth
\lower1.5ex\hbox{`}\hidewidth\crcr\unhbox0}}} \def\cprime{$'$}
\def\cprime{$'$} \def\cprime{$'$}
\def\polhk#1{\setbox0=\hbox{#1}{\ooalign{\hidewidth
\lower1.5ex\hbox{`}\hidewidth\crcr\unhbox0}}} \def\cprime{$'$}
\def\cprime{$'$} \def\polhk#1{\setbox0=\hbox{#1}{\ooalign{\hidewidth
\lower1.5ex\hbox{`}\hidewidth\crcr\unhbox0}}} \def\cprime{$'$}
\def\cprime{$'$} \def\cydot{\leavevmode\raise.4ex\hbox{.}} \def\cprime{$'$}
\def\cprime{$'$} \def\cprime{$'$} \def\cprime{$'$}

\end{document}